\documentclass{amsart}

\usepackage{geometry}
\usepackage{graphicx}
\usepackage{amsmath,amsthm,amssymb}
\usepackage{float}
\usepackage{hyperref}
\hypersetup{
    colorlinks=true,
}

\newcommand{\R}{\mathbb{R}}
\newcommand{\iim}{\mathrm{i}}

\newtheorem{thm}{Theorem} 
\newtheorem{lem}[thm]{Lemma} 
\newtheorem{rmk}{Remark} 
\newtheorem{prop}[thm]{Proposition}

\begin{document}

\title{Convergence analysis of a Schr\"{o}dinger problem with moving boundary}

\author{Daniel G. {Alfaro~Vigo}}%\footnote{Corresponding author}}
\address{Institute of Computing \& Graduate Program in Informatics, Federal University of Rio de Janeiro, Caixa Postal 68530, CEP 21941-590, Rio de Janeiro - RJ, Brazil}
%\thanks{Corresponding author}
\email{dgalfaro@ic.ufrj.br}
\author{Daniele C.R. Gomes}
\address{Graduate Program in Informatics, Federal University of Rio de Janeiro, Caixa Postal 68530, CEP 21941-590, Rio de Janeiro - RJ, Brazil}
\email{daniele.rocha@ppgi.ufrj.br}
\author{Bruno A. do Carmo}
\address{Graduate Program in Informatics, Federal University of Rio de Janeiro, Caixa Postal 68530, CEP 21941-590, Rio de Janeiro - RJ, Brazil}
\email{bruno.carmo@ppgi.ufrj.br}
\author{Mauro A. Rincon}
\address{Institute of Computing \& Graduate Program in Informatics, Federal University of Rio de Janeiro, Caixa Postal 68530, CEP 21941-590, Rio de Janeiro - RJ, Brazil}
\email{rincon@ic.ufrj.br}


\date{\today}

\begin{abstract}
In this article, we present the mathematical analysis of the convergence of the linearized Crank-Nicolson Galerkin method for a nonlinear Schrödinger problem related to a domain with a moving boundary. 
The convergence analysis of the numerical method is carried out for both semidiscrete and fully discrete problems. 
We establish an optimal error estimate in the $L^2$-norm with order $\mathcal{O}(\tau^2+ h^s),~ 2\leq s\leq r$, where $h$ is the finite element mesh size, $\tau$ is the time step, and $r-1$ represents the degree of the finite element polynomial basis.
Numerical simulations are provided to confirm the consistency between theoretical and numerical results, validating the method and the order of convergence for different degrees $p\geq 1$ of the Lagrange polynomials and also for Hermite polynomials (degree $p=3$), which form the basis of the approximate solution. 
\end{abstract}

\maketitle

%% main text
\section{Introduction}\label{sec:intro}

Let $T>0$ and $\Omega$ be an open bounded set of $\R^n$, with a regular boundary $\Gamma$ of class $C^2$. 
We represent by $k(t)$ a real function defined on the set of all non-negative real numbers $[0,\infty)$.
Consider the subset $\Omega_{t}=\{ x\in\R^n : x=k(t)y,\, y \in \Omega\}\subset\R^n$, for each $t\in[0,T]$. Let $\widehat{Q}$ be the noncylindrical domain of $\R^{n+1}$ with regular one-sided boundary $\widehat{\Sigma}$ defined by
\begin{align*}
\widehat{Q}={\bigcup\limits_{0<t<T}}\{\Omega_{t}\times\{t\}\},\quad 
\widehat{\Sigma}={\bigcup\limits_{0<t<T}}\{\Gamma_{t}\times\{t\}\}, \quad\mbox{where}\quad \Gamma_{t}=\partial\Omega_{t}.
\end{align*}

We shall consider the following nonlinear Schr\"{o}dinger problem in a noncylindrical domain:
\begin{equation}
\label{prob:principal}
\quad
\left\{ %%%%% 
\begin{aligned}
& \partial_{t} u(x,t) - \iim \Delta u(x,t) + \left|u(x,t)\right|^{\rho} u(x,t) = \hat{f}(x,t) \ \ \mbox{in} \ \ \hat{Q},
\\
& u(x,t) = 0 \ \ \mbox{on} \ \ \hat{\Sigma},
\\
& u(x,0)=u_0(x) \ \ \mbox{in} \ \ \Omega_0,
\end{aligned}
\right.
\end{equation}
%%%%%%
where $\partial_{t} w$ denotes the partial derivative of function $w$ with respect to time $t$, $\iim$ represents the imaginary unit and $\Delta$ represents the Laplacian operator. 

Let us consider the diffeomorphism $\mathcal{T}:\widehat{Q}\rightarrow Q$ 
defined by $(y,t)=\mathcal{T}(x,t)=(x/k(t), t)$, see, for example, \cite{Artigo2003Raquel,Artigo2003Liu}.
By considering the function $v(y,t)=u(k(t)y,t)$, we transform the problem \eqref{prob:principal} into the following equivalent one
\begin{equation} \label{prob:equivalente}
\quad\left\{
\begin{aligned}
& \partial_{t} v(y,t) - \frac{k'(t)}{k(t)} y_j \frac{\partial v(y,t)}{\partial y_j} -  \frac{\iim}{k^2(t)}\Delta v(y,t) 
+|v(y,t)|^{\rho}v(y,t) = f(y,t) ~\mbox{in}~ Q, 
\\
& v(y,t) = 0 \ \ \mbox{on} \ \ \Sigma, 
\\
& v(y,0)=v_0 \ \ \mbox{in} \ \ \Omega, 
\end{aligned}
\right.
\end{equation} 
where $f(y,t)=\hat{f}(k(t)y,t)$ and when $(x,t)$ varies in $\widehat{Q}$ the point $(y,t)$, with $y=x/k(t)$, varies in $Q=\Omega\times(0,T)$. 
%%%%%%
Thus, we obtain a nonlinear Schr\"{o}dinger equation with advection and time-dependent coefficients in a cylindrical domain.

%%%%%%%%%%%

The nonlinear Schrödinger equation is widely used to model a broad spectrum of physical phenomena, as it describes the envelope dynamics of a quasi-monochromatic plane wave propagating in a weakly nonlinear dispersive medium. Its applications span various fields, including superconductivity, nonlinear optics, Bose-Einstein condensation, and surface gravity waves. See, for example, \cite{Dalfovo1999, Hartree, Kelley, Cipolatti, Strecker2003}. 
Furthermore, the externally forced nonlinear Schrödinger equation arises in many physical situations such as charge density waves, optical fibers, and plasmas driven by rf fields (see, for instance, \cite{Kaup_1978, Malomed_1995, Nozaki_1986}.). 

More recently, the study of the nonlinear Schrödinger equation with a moving boundary (i.e., in a noncylindrical domain) has garnered significant attention from both mathematical and physical perspectives. From a mathematical point of view, the existence and uniqueness of solutions were addressed in \cite{Artigo2019Gomes, Artigo2024Gomes}, while \cite{Bisognin_etal} focused on the asymptotic properties of the solutions and \cite{Beauchard_etal, Castelli2016} explored their controllability properties. From a physical perspective, studies have focused mainly on Bose–Einstein condensates (see, for instance, \cite{Band_etal, Theodorakis2009, Campo2012, VanGorder2021}).  

It is worth pointing out that several important applications that require the control of Bose-Einstein condensates have been proposed, including quantum interferometry, metrology, and quantum computation (see, for example, \cite{Pinsker2015, Ray_etal, Sabin_2014, Frank_etal}). Moreover, the aforementioned results suggest that this control can be achieved through the use of a moving boundary, further emphasizing the importance of studying the nonlinear Schrödinger equation in noncylindrical domains.

Numerical methods for the nonlinear Schrödinger equation in cylindrical domains using finite difference, finite element, spectral representations, etc., have been extensively applied and theoretically studied. See, for example, \cite{Akrivis1991, Yves1991, Chang_etal, Jilu2013, Shi2016, X.Li_etal, Zhang2017, Asadzadeh2019, Li_etal2021, Ibarra-Villalon_etal, Yang2023} and references therein.

However, numerical studies on the nonlinear Schrödinger equation with a moving boundary are relatively rare. In \cite{Artigo2019Gomes, Artigo2024Gomes}, the authors presented theoretical results on the existence and uniqueness of solutions to \eqref{prob:principal} and conducted numerical simulations in one- and two-dimensional settings using a linearized Crank-Nicolson Galerkin method. Their numerical experiments demonstrated an optimal convergence order of the approximate solutions, but no theoretical analysis of the numerical method was provided. 
The recent preprint \cite{mollisaca_etal} presents theoretical results for a higher-order nonlinear Schrödinger equation in a one-dimensional moving domain, along with numerical simulations performed using a finite difference method and an analysis of its properties. %
In \cite{VanGorder2021}, the author performed numerical simulations to illustrate their theoretical findings, focusing on different methods for controlling Bose-Einstein condensates.
Reference \cite{Ljungberg_etal} presents a finite difference scheme for the linear Schrödinger equation in a moving domain.

We note that to develop their numerical studies, these references applied a time-dependent coordinate transformation, resulting in a linear or nonlinear Schrödinger equation with advection and time-dependent coefficients (see, for example, problem \eqref{prob:equivalente}).

The main objective of this paper is to present a detailed convergence analysis of the linearized Crank-Nicolson Galerkin method presented in \cite{Artigo2019Gomes, Artigo2024Gomes} for the nonlinear Schrödinger equation \eqref{prob:principal} in the case where $n = 1, 2, 3$. The Crank-Nicolson method and its variants are commonly used to efficiently and accurately solve the nonlinear Schrödinger equation, (see e.g. \cite{Akrivis1991, Chang_etal, Shi2016,  Zhang2017, Li_etal2021, Yang2023, Li_etal2025}), and other types of partial differential equations (e.g. \cite{Artigo1973_Wheeler, Livro2006_Thomee,  ArtigoNatanael2016, Zhang2019}).

The major contributions of this paper are discussed below.
\begin{itemize}
\item We derive optimal error estimates of the approximate solutions to the problem \eqref{prob:principal}. More specifically, for the linearized Crank-Nicolson Galerkin method  we obtain an error estimate of order $\mathcal{O}(\tau^2 + h^r)$ in the $L^2$-norm, where $h$ represents the finite element mesh size, $r - 1$ the degree of the polynomial basis and $\tau$ the time step, under the condition $\tau = {o}(h^{n/4})$. 
Additionally, for the semidiscrete Galerkin finite element method we derive error estimates of the orders $\mathcal{O}(h^r)$ and $\mathcal{O}(h^{r-1})$ in the $L^2$ and $H^1$ norms, respectively.

To the best of our knowledge, such estimates for a nonlinear Schrödinger equation with advection and time-dependent coefficients have not been reported before. For the generalized nonlinear Schrödinger equation with a nonlinearity of the form $f(|u|^2)u$, where $f \in C^2(\R)$, optimal error estimates in the $H^1$ norm for the Crank-Nicolson Galerkin method, under the time step condition $\tau = o(h^{n/4})$, were presented in \cite{Yves1991}. Similarly, optimal error estimates in the $L^2$ norm were obtained for the linearized Crank-Nicolson Galerkin method with nonconforming finite elements under the time step condition $\tau = \mathcal{O}(h)$ in \cite{Zhang2017}, as well as for the Crank-Nicolson Galerkin method without any time step restriction in \cite{Jilu2013}.  For the Schrödinger equation with cubic nonlinearity optimal error estimates without any time step restriction were presented in \cite{Yang2023} for a modified Crank-Nicolson Galerkin method. 
\item We introduce a new time-dependent Ritz projection, which plays a critical role in establishing the optimal convergence order of the method. This projection allows for a straightforward treatment of the advection term and the time-dependent coefficients in \eqref{prob:equivalente}.
The Ritz projection is extensively used to establish error estimates of finite element methods. 
\item We establish uniform bounds for the approximate solutions in the $L^\infty$-norm, which, combined with an induction argument, are crucial for handling the nonlinear term. In contrast, other studies rely on a priori bounds for the nonlinear term that depend on the discretization parameters or assume the validity of appropriate bounds of approximate solutions; for example in \cite{Artigo1973_Wheeler, Artigo2013Kirby, ArtigoNatanael2016, Shi2016, Asadzadeh2019}.

\end{itemize}

This article is organized as follows. In Section \ref{sec:exist}, we briefly describe theoretical results on the existence and uniqueness of solutions to problem \eqref{prob:principal}. In Section \ref{sec:num_method}, we introduce the linearized Crank–Nicolson Galerkin method used to obtain approximate numerical solutions to this problem. In Section \ref{sec:num_an},  error estimates for semidiscrete and fully discrete problems are presented. In Section \ref{sec:num_ex}, we present some numerical examples to illustrate the validity of the theoretical analysis. Finally, in Section \ref{sec:conc}, we present some final remarks.

\section{Existence and uniqueness of solutions}\label{sec:exist}

We formulate a recent result concerning the existence and uniqueness of solutions to problem \eqref{prob:principal}. For this analysis, we introduce the following hypotheses:\\
\begin{enumerate}
\item[{\bf H1:}] $k \in W_{loc}^{2,\infty}([0,\infty[); k(t) \geq k_0 > 0, \ \forall t \geq 0$;
\item[{\bf H2:}] $0\leq\rho<\infty$ ~if $n=1,2$ and $0\leq\rho\leq\dfrac{2}{n-2}$ if $n\geq 3$. \\
\end{enumerate}

\begin{thm}
\label{teo:exi.uni.u}
We assume that hypotheses (H1) and (H2) are satisfied. 
Let us consider the initial data $u_{0}\in H_{0}^{1}(\Omega_{0})$ and $\widehat{f}\in L^{2}(0,T;H_{0}^{1}(\Omega_{t}))$; then
there exists a function $u:\widehat{Q}\longrightarrow\mathbb{C}$ such that
\begin{enumerate}
\item
$u \in L^\infty(0,T;H_0^1(\Omega_t))\cap L^p(0,T;L^p(\Omega_t)), \ \mbox{with} \ p=\rho+2;$
%%%%%%
\item
$\partial_{t} u \in L^{p'}(0,T;H^{-1}(\Omega_t)), \ \mbox{with} \ p'=\frac{\rho+2}{\rho+1};$
\item
$\partial_{t} u - \iim \Delta u + |u|^{\rho} u = \hat{f} \ \mbox{in} \ L^{p'}(0,T;H^{-1}(\Omega_t));$
\item
$u(0)=u_0 \ \mbox{in} \ \Omega_0.$
\end{enumerate} 
\end{thm}

Due to diffeomorphism $\mathcal{T}$, we know that $u$ is a solution to problem \eqref{prob:principal} given by Theorem \eqref{teo:exi.uni.u} if and only if $v$ is a solution to problem \eqref{prob:equivalente} given by the following theorem:

\begin{thm}
\label{teo:exi.uni.v}
We assume that hypotheses (H1) and (H2) are satisfied. 
Let us consider the initial data $v_0 \in H_0^1(\Omega)$ and $f \in L^2(0,T;H_0^1(\Omega))$, then
there exists a function $v:Q \rightarrow \mathbb{C}$, such that
\begin{enumerate}
\item
$v \in L^\infty(0,T;H_0^1(\Omega))\cap L^p(0,T;L^p(\Omega)), \ \mbox{with} \ p=\rho+2;$ 
%%%%%%
\item
$\partial_{t} v \in L^{p'}(0,T;H^{-1}(\Omega)), \ \mbox{with} \ p'=\frac{\rho+2}{\rho+1};$
\item
$\partial_{t} v - \dfrac{k'}{k} y_j \dfrac{\partial v}{\partial y_j} -  \dfrac{\iim}{k^2}\Delta v + |v|^{\rho}v = f \ \mbox{in} \ L^{p'}(0,T;H^{-1}(\Omega));$
\item
$v(0)=v_0 \ \mbox{in} \ \Omega.$
\end{enumerate} 
\end{thm}

The above equation can be written as 
%%%%%%
\[
\partial_{t} v + \frac{\iim}{k^2} {L} v +  \left|v\right|^{\rho}v = f %
\]
where we considered the operator
\[
{L} u = - \Delta u + \iim \displaystyle (k'k) y \cdot \nabla u. 
\]

\begin{rmk}\label{rem:num_adv}
The numerical solution to problem \eqref{prob:principal} can be readily obtained from problem \eqref{prob:equivalente} which  is defined in a cylindrical domain. This greatly simplifies the development of a numerical method, as there is no need to include a remeshing process.   
\end{rmk}

%%%%%%%%%%%%%%%
%%%%%%%%%%%%%%%

\section{Numerical methods}\label{sec:num_method}

In order to obtain approximate solutions to \eqref{prob:principal}, we apply a linearized Crank-Nicolson Galerkin method to problem \eqref{prob:equivalente}. This method is based on a finite element Galerkin discretization of the spatial variable and a modified Crank-Nicolson time-stepping scheme.

%------------------------------------------------------------------------
\subsection{Semidiscrete problem}
%------------------------------------------------------------------------

The starting point for the application of the finite element Galerkin approximation is the variational formulation of problem \eqref{prob:equivalente}. % 
We introduce the time-dependent sesquilinear form on $H_0^1(\Omega)$, associated with operator $L$, given as
\begin{equation}\label{eq:sesq_lin}
a_L(t;u,w) =   (\nabla u,\nabla w)  + \iim \gamma(t)(y \cdot \nabla u,w)
\end{equation}
where $\gamma(t) = k (t) k'(t)$. 
Then, the variational formulation consists on finding $v:[0,T]\rightarrow H_0^1(\Omega)$, which satisfies the initial condition $v(0)=v_{0}$  and such that
\begin{equation}
\label{prob:variational}
%%%%%
(\partial_{t} v(t),w) + \frac{\iim }{k^2(t)}a_L(t; v(t), w) + (g(v(t)),w) = (f(t),w),
\end{equation}
 for all $w\in H_0^1(\Omega)$ and  $t\in[0,T]$. 
 In order to simplify the notation, we introduced the function
$g:\mathbb{C}\rightarrow \mathbb{C}, z \mapsto g(z)=|z|^{\rho}z$.

Let $\mathcal{T}_h$ be a quasiuniform finite element discretization of  $\Omega_h\subset\Omega$ consisting of the finite union of the elements $K\in\mathcal{T}_h$. The parameter $h$ denotes the maximal diameter of the elements of $\mathcal{T}_h$. %$\Omega_h$ represents a polygonal domain  
Additionally, there is a finite dimensional function space $V_h = S_h^r(\Omega)\subset H^1_0(\Omega)$ of continuous function on $\overline{\Omega}$ and for each $K\in\mathcal{T}_h$ the subspace $P_K = \{v_{h|K}, v_h\in V_h\}$ consists of polynomials of degree less than or equal $r-1$ ($r\geq2$). In the space $S_h^r(\Omega)$, we can consider the norm  $\|\cdot\|_{H^1_0(\Omega)}$, which is also equivalent to the norms $\|\cdot\|_{L^2(\Omega)}$ and $\|\cdot\|_{L^\infty(\Omega)}$ with constants that can depend on $h$.

The solution to problem \eqref{prob:variational} is approximated by a function $v_h:[0,T]\rightarrow S_h^r(\Omega)$, such that for all $\chi\in S_h^r(\Omega)$ and $t\in[0,T]$
\begin{equation}
\label{prob:aproximado}
%%%%%%
(\partial_t v_h(t),\chi) + \frac{\iim }{k^2(t)}a_L(t; v_h(t), \chi) + (g(v_h(t)),\chi) = (f(t),\chi),
\end{equation}
with the initial condition $v_h(0)=v_{0h}$ where  $v_{0h}\in S_h^r(\Omega)$  is an approximation of $v_{0}$. 

\begin{rmk}
The semidiscrete problem \eqref{prob:aproximado} can be recast as a nonlinear system of first order ordinary differential equations for the coefficients of the expansion of $v_h$ with respect to a  basis of $S_h^r(\Omega)$. The existence and uniqueness of a solution can be obtained from Caratheodory's theorem and the local Lipschitz continuity of the right-hand side of the system.
\end{rmk}

%%%%%%%%%%%%

%------------------------------------------------------------------------
\subsection{Fully discrete problem}
%------------------------------------------------------------------------
Since the general solution to problem \eqref{prob:aproximado} can not be obtained explicitly, we consider a time-stepping method based on the Crank-Nicolson scheme.  This method relies on a linearization that avoids the appearance of nonlinear equations in the time-stepping procedure. We refer to this method as a linearized Crank-Nicolson Galerkin scheme.

Let $N\in\mathbb{N}$ and take a uniform discretization $0=t_0<t_1<\cdots<t_N=T$ of the interval $[0,T]$ with a time step of size $\tau =T/N$.
 
Let $U^m$ be an approximation of $v(t_m)$ in $S_h^r(\Omega)$, where $t_m=m\tau$ ($0\leq m \leq N$).
The fully discrete problem consists in finding $\{U^m\}_{m=0}^N$ from $S_h^r(\Omega)$, such that $m\in\{2, \ldots, N\}$,
\begin{align}\label{prob_aprox:discrete}
 %%%%%
 (\delta_\tau U^m, \chi) 
+ \frac{\iim}{k^2_{m-\frac{1}{2}}} a_L(t_{m-\frac{1}{2}}; \widehat{U}^m, \chi)
+ (g(\widetilde{U}^m),\chi)
= (f_{m-\frac{1}{2}}, \chi),
\end{align}
for all $\chi \in S_h^r(\Omega)$, where $t_{m-\frac{1}{2}}=(t_m +t_{m-1})/2$ is the midpoint of interval $[t_{m-1}, t_{m}]$ and $k^2_{m-\frac{1}{2}} = k^2(t_{m-\frac{1}{2}})$, $f_{m-\frac{1}{2}}=f(t_{m-\frac{1}{2}})$
\begin{align}\label{eq:discrete_diff}
%%%%%
\delta_\tau U^{m} = \dfrac{U^m-U^{m-1}}{\tau}, \quad
\widehat{U}^{\; m} = \dfrac{U^m+U^{m-1}}{2}, \quad
\widetilde{U}^{\; m} = \dfrac{3U^{m-1}-U^{m-2}}{2}. 
\end{align}
We notice that the initial approximations $U^0$ and $U^1$ shall be given; their choice will be discussed later.

It is worth noticing how the linearized Crank-Nicolson Galerkin scheme can be established. The variational problem \eqref{prob:variational} at time $t=t_{m-\frac{1}{2}}$, reads
\begin{align*}
 %%%%%
(\partial_{t} v(t_{m-\frac{1}{2}}),w) 
+ \frac{\iim}{k^2_{m-\frac{1}{2}}}
   a_L(t_{m-\frac{1}{2}}; v(t_{m-\frac{1}{2}}), w) 
+ (g(v(t_{m-\frac{1}{2}})),w)
= (f_{m-\frac{1}{2}},w), 
\end{align*}
for all $ w\in H_0^1(\Omega)$. Using Taylor expansion, we have the following second order approximations %at time $t_{n-\frac{1}{2}}$ 
\begin{align*}
 %%%%%
 \partial_{t} v(t_{m-\frac{1}{2}})
= \dfrac{v(t_{m})-v(t_{m-1})}{\tau} 
+ \mathcal{O}(\tau^2),
 \quad 
  v(t_{m-\frac{1}{2}})
= \dfrac{v(t_{m})+v(t_{m-1})}{2} 
+ \mathcal{O}(\tau^2).
\end{align*}
This leads to the problem of finding $\{U^m\}_{m=0}^N$ from $S_h^r(\Omega)$, such that for $m\in\{1, \ldots, N\}$,
\begin{align*}
 %%%%%
 (\delta_\tau U^m, \chi) 
+ \frac{\iim}{k^2_{m-\frac{1}{2}}}
  a_L(t_{m-\frac{1}{2}}; \widehat{U}^m, \chi)
+ (g(\widehat{U}^m),\chi)
= (f_{m-\frac{1}{2}}, \chi), 
\end{align*}
for all $\chi \in S_h^r(\Omega)$. This is the so-called Crank-Nicolson Galerkin method. 
Notice that in this method, we need to solve a nonlinear system of equations at each time step.

In order to avoid this situation, the following approximation is used for the nonlinear term
\begin{align*}
  v(t_{m-\frac{1}{2}})
= \frac{3v(t_{m-1})-v(t_{m-2})}{2} 
+ \mathcal{O}(\tau^2).
\end{align*}
This modification leads to the linearized Crank-Nicolson Galerkin scheme given by \eqref{prob_aprox:discrete}. In this method we obtain a linear system of equations at each time step. However, this method requires prior knowledge of the approximation  $U^1$.

To obtain $U^1$, we choose a predictor-corrector one-step method discussed in \cite{Livro2006_Thomee, Artigo1973_Wheeler} in the context of nonlinear parabolic problems. The method is defined as follows. 

Let $U^{1^-}\in S_h^r(\Omega)$ be the predictor for $v(t_1)$.
Taking $m=1$ in equation \eqref{prob_aprox:discrete} and substituting $\widetilde{U}^1$ by $U^0$, that is approximating $v(t_{\frac{1}{2}})$ by $v^0$ in the nonlinear term, we get 
\begin{align}\label{prob_aprox:U^0(1)}
  \Big(\dfrac{U^{1^-}-U^0}{\tau}, \chi\Big) 
+ \frac{\iim}{k^2_{\frac{1}{2}}}
  a_L\Big(t_{\frac{1}{2}}; \dfrac{U^{1^-}+U^0}{2}, \chi\Big) 
+ (g({U}^0),\chi)
= (f_{\frac{1}{2}}, \chi), 
\end{align}
for all $\chi \in S_h^r(\Omega)$.

The corrector $U^{1}\in S_h^r(\Omega)$ is obtained using \eqref{prob_aprox:discrete} for $m=1$ and substituting $\widetilde{U}^1$ by $\dfrac{U^{1^-}+U^0}{2}$, i.e. approximating $v(t_{\frac{1}{2}})$ by $\dfrac{U^{1^-}+U^0}{2}$ in the nonlinear term. We get that
\begin{equation}
\label{prob_aprox:U^0(2)}
 %%%%%
 (\delta_\tau U^1, \chi)
+ \frac{\iim}{k^2_{\frac{1}{2}}}
  a_L(t_{\frac{1}{2}}; \widehat{U}^1,\chi)
+ \Big(g\Big(\dfrac{U^{1^-}+U^0}{2}\Big),\chi\Big) 
= (f_{\frac{1}{2}},\chi), 
\end{equation}
for all $\chi \in S_h^r(\Omega)$.

The existence and uniqueness of the solutions to problems \eqref{prob_aprox:discrete}, \eqref{prob_aprox:U^0(1)} and \eqref{prob_aprox:U^0(2)} is guaranteed for any sufficiently small $\tau$.

We can rewrite problems \eqref{prob_aprox:discrete}, \eqref{prob_aprox:U^0(1)} and \eqref{prob_aprox:U^0(2)} in the following unified form. Find $U^m\in S_h^r(\Omega)$ for $m=1^-,1,2, \dots,N$ such that 
\begin{equation}
\label{prob_aprox:discrete_all}
%%%%%
(\delta_\tau U^m, \chi) +  \frac{\iim}{k^2_{m-\frac{1}{2}}} a_L(t_{m-\frac{1}{2}}; \widehat{U}^m, \chi)+(g^m,\chi)= (f_{m-\frac{1}{2}}, \chi),
\end{equation}
for all $\chi \in S_h^r(\Omega)$, where we introduced the notation 
\begin{equation}\label{eq:notation_all}
%%%%%%
\delta_{\tau}U^{1^-} = \frac{U^{1^-}-U^{0}}{\tau},\quad %\\ 
\hat{U}^{1^-} = \frac{U^{1^-}+U^{0}}{2},\quad %\\
g^m = 
\begin{cases}
g(U^0),&\text{if $m=1^-$}\\
g(\hat{U}^{1^-}),&\text{if $m=1$}\\
g(\tilde{U}^m),&\text{if $m\geq 2$} %
\end{cases}
\end{equation}
and also considered that $1^- -\frac{1}{2} = \frac{1}{2}$.
%

%%%%%%%%%%%%%
%%%%%%%%%%%%%
\section{Numerical analysis}\label{sec:num_an}

As mentioned before, we shall consider two cases for the error analysis. In the first case, we approximate the function $v(y, t)$ by the function $v_h(y, t)$, which, for each fixed $t$, belongs to the approximate finite-dimensional space $S_h^r(\Omega)$ defined in section \ref{sec:num_method}. This function will be a solution to a  Galerkin approximate problem that consists of a system of ordinary differential equations in time, referred to as the semidiscrete problem. In other words, since time varies continuously, the error analysis for the semidiscrete problem will be in space. For this analysis, we will use the variational formulation \eqref{prob:variational}. In the second part, we analyze the fully discrete problem.

%%%%%%%%%%%%%%%%%%%
\subsection{Auxiliary results}
Before we proceed to the presentation of error estimates, we introduce some auxiliary results.

\begin{prop} \label{prop:sesq_lin}
The sesquilinear form \eqref{eq:sesq_lin} is continuous  in $H^1_0(\Omega)$ and coercive in $H^1_0(\Omega)$ with respect to $ L^2(\Omega)$ uniformly in $t$, i.e. there are constants $C>0$, $\lambda_0 >0$ such that for  $t\in[0,T]$ and $u, w\in H^1_0(\Omega)$
\begin{align}\label{eq:coercive}
&|a_L(t; u, w)| \leq C\|u\|_{H^1_0(\Omega)} \|w\|_{H^1_0(\Omega)},  \\
&Re\,a_L(t; u, u) + \lambda_0 \|u\|_{L^2(\Omega)} ^2\geq  \tfrac{1}{2} \|u\|_{H^1_0(\Omega)}^2.
\end{align}
\end{prop}

The proof of this Proposition is straightforward and will be omitted.

\begin{rmk}
We can take $C=\sqrt{1+ d_\Omega \|kk'\|_\infty}$ and $\lambda_0 = \tfrac{(d_\Omega \|kk'\|_\infty)^2}{2}$ where $d_\Omega = \sup\{|y|: y\in\Omega\}$.
\end{rmk}

As a direct consequence of Proposition \ref{prop:sesq_lin}, we get that the sesquilinear  form 
\[
a_{L_0}(t; u, w) = a_{L}(t; u, w) + \lambda_0 (u, w),
\] 
is elliptic and continuous on $H^1_0(\Omega)$ for each fixed $t\in[0,T]$. This form is associated with the elliptic operator $L_0 u = L  u + \lambda_0 u$.

We define a time-dependent Ritz projection $P_h$ from $H^1_0(\Omega)$ onto $S_h^r(\Omega)$ ($P_h:[0,T]\times H^1_0(\Omega)\to S_h^r(\Omega)$) such that for $t\in[0,T]$ and $u\in H^1_0(\Omega)$
\begin{equation}\label{eq:t-ritz_proj}
a_{L_0}(t; P_h(t, u), \chi) = a_{L_0}(t; u, \chi),\quad\forall \chi\in S_h^r(\Omega).
\end{equation}
The existence of this projection is guaranteed by the Lax-Milgram theorem since $a_{L_0}$ is elliptic and continuous on $H^1_0(\Omega)$ for each fixed $t\in[0,T]$. Similarly, we define the adjoint Ritz projection $P_h^*:[0,T]\times H^1_0(\Omega)\to S_h^r(\Omega)$ such that
\[
a_{L_0}(t; \chi, P_h^*(t, u)) = a_{L_0}(t; \chi, u),\quad\forall \chi\in S_h^r(\Omega).
\]
Additionally, for any function $u:[0,T]\to H^1_0(\Omega)$ we introduce the function $P_hu: [0,T]\to  S_h^r(\Omega)$ such that $t\mapsto P_h(t, u(t))$. % 

In the case where $k(t) = \mathrm{const}$, $P_h$ coincides with the ``classical" Ritz projection $R_h$ which solves the variational problem $(\nabla R_h u ,\nabla\chi) = (\nabla u ,\nabla\chi)$ for all $\chi\in S_h^r(\Omega)$ (see \cite{Livro2002_Ciarlet}). The projection $R_h$ satisfies several properties widely used to obtain error estimates in parabolic and elliptic problems (see, for instance, \cite{Livro2002_Ciarlet, Livro2006_Thomee}).   In general, when $k(t) \neq \mathrm{const}$,  projection $P_h$ satisfies similar properties. 

\begin{lem}\label{lem:ritz_proj}
The Ritz projection $P_h$ has the following properties
\begin{enumerate}\renewcommand{\theenumi}{\roman{enumi}}%
\item Let $u\in H^1_0(\Omega)\cap H^s(\Omega)$ with $1\leq s\leq r$ then there exists a constant $C>0$ such that for $t\in[0,T]$
\begin{equation}\label{eq:est_proj_i}
\|u-P_h(t, u)\|_{L^2(\Omega)} + h\|u-P_h(t, u)\|_{H_0^1(\Omega)} \leq C h^s \|u\|_{H^s(\Omega)}.
\end{equation}
\item Let $u\in H^1_0(\Omega)\cap H^s(\Omega)$ where $1\leq s\leq r$ for $n=1$ and $2\leq s\leq r$ for $n=2,3$. Then, there is a constant $C>0$ such that for $t\in[0,T]$
\begin{equation}\label{eq:est_proj_ii}
\|u-P_h(t, u)\|_{L^\infty(\Omega)}  \leq C \ell_h h^{s-n/2} \|u\|_{H^s(\Omega)}% \|u\|_{H^s(\Omega)}
\end{equation}
where 
\[
\ell_h = 
\begin{cases}
1 + |\log h|,& \text{if $r=2$},\\
1, & \text{if $r>2$}.
\end{cases}
\]
\item  Let 
$u, \partial_{t} u,\dots,\partial_{t}^{m} u\in L^2(0,T;H^1_0(\Omega)) $ 
and $k\in C^{m+1}([0,T])$ with $m\geq 1$, then for $j=1,\dots,m$ the derivatives $z_j(t) = (P_h u)^{(j)}(t) = \partial_{t}^j (P_h u)(t)$ exist for almost every  $t\in(0,T)$ and satisfy the equation 
\begin{align}\label{eq:def_DP}
& a_{L_0}(t; z_{j}(t), \chi)
= a_{L_0}(t; \partial_t^j u(t), \chi) 
+ i\sum_{l=0}^{j-1} \dbinom{j}{l} \gamma^{(j-l)}(t) (y\cdot\nabla [\partial_t^l u(t)-z_{l}(t)],\chi), 
\end{align}
 for all $\chi\in S_h^r(\Omega)$, where $\gamma(t) = k (t) k'(t)$ and $\binom{j}{l}=\frac{j!}{l! (j-l)!}$. 
\item Let 
$u, \partial_{t} u,\dots,\partial_{t}^{m} u\in L^2(0,T;H^1_0(\Omega)\cap H^s(\Omega))$, $k\in C^{m+1}([0,T])$ with $m\geq 1$ and $1\leq{s}\leq{r}$. Then, there exists a  constant $C>0$ such that for $j=1,\dots,m$ and almost every $t$
\begin{align}\label{eq:est_proj_iv}
& \|\partial_t^j u(t)-\partial_t^j (P_h u)(t)\|_{L^2(\Omega)} 
+ h\|\partial_t^j u(t)-\partial_t^j (P_h u)(t) \|_{H_0^1(\Omega)}
\leq 
C h^s  \sum_{j=0}^m\|\partial_t^j u(t)\|_{H^s(\Omega)}.
\end{align}
\end{enumerate}
\end{lem}

\begin{rmk} 
Evidently,  the adjoint Ritz projection $P_h^*$ also has the properties presented in Lemma \ref{lem:ritz_proj}.
It is also worth noticing that in general,
$\partial_t^j (P_h u)(t) \neq P_h(t; \partial_t^j u(t))$; however, we will establish that they are asymptotically close (see equation \eqref{eq:ineq_2}).
\end{rmk}

\begin{rmk}\label{obs:prop_iii-iv}
Properties (iii) and (iv) deal with functions lying in spaces of the form 
\[
W^{m}(0,T;X) =\bigl\{u: (0,T)\to X\;\bigl|\; 
\partial_t^j u \in L^2(0,T;X),\quad 0\leq j\leq m\bigr\}
\]
where $X$ is a separable Hilbert space and a norm in $W^{m}(0,T;X)$ is given as  $$\|u\|_{W^{m}(0,T;X)} = \left[\sum_{j=0}^m \|
\partial_t^j u \|^2_{L^2(0,T;X)}\right]^{1/2}.$$ The space $W^{m}(0,T;X)$ is continuously imbedded in $C^{m-1}([0,T];X)$. Moreover, if $u\in W^{m}(0,T;X)$, then 
$\partial_t^{m-1} u(t)$ is absolutely continuous and differentiable almost everywhere (see, for instance, \cite{Cazenave2003}).
\end{rmk}

\begin{proof}%[of lemma \ref{lem:ritz_proj}] 
The proof of these properties is based on some ideas used to establish Céa and Aubin-Nitsche lemmas (see, for instance, \cite{Livro2002_Ciarlet, Livro2006_Thomee}).

(i)\quad Let $t\in[0,T]$, $u\in H^1_0(\Omega)\cap H^s(\Omega)$, set $u_h = P_h(t; u)\in S_h^r(\Omega)$ and take $\chi\in S_h^r(\Omega)$.  From the definition of the Ritz projection, equation \eqref{eq:coercive} and the properties of the space $S_h^r(\Omega)$,  we get that
\begin{align*}
 \|u-u_h\|_{H^1_0(\Omega)}^2 &\leq  2 Re\, a_{L_0}(t; u-u_h, u-u_h) \\
&\leq 2 Re\,\left\{ a_{L_0}(t; u-u_h, u-\chi) + a_{L_0}(t; u-u_h, \chi-u_h) \right\}\\
&\leq 2 Re\, a_{L_0}(t; u-u_h, u-\chi) \\
&\leq C \|u-u_h\|_{H^1_0(\Omega)} \|u-\chi\|_{H^1_0(\Omega)}.
\end{align*}
Consequently,
\begin{equation}\label{eq:ineq_1}
 \|u-u_h\|_{H^1_0(\Omega)} \leq C\inf_{\chi\in H^1_0(\Omega)} \|u-\chi\|_{H^1_0(\Omega)} \leq C_1 h^{s-1} \|u\|_{H^s(\Omega)}.
\end{equation}

Consider the operator $L_0^*$, adjoint to $L_0$. Let $\varphi = u-u_h$ and take $\psi$ that solves the equation $L_0^* \psi = \varphi$ with the boundary condition $\psi|_{\partial\Omega} = 0$. Since $L_0^*$ is a second order elliptic operator\footnote{We assume that the boundary $\partial\Omega$ is sufficiently smooth.}, $\psi\in H^2(\Omega)\cap H_0^1(\Omega)$ and $\|\psi\|_{H^2(\Omega)} \leq C_1\|\varphi\|_{L^2(\Omega)}$ for some constant $C_1>0$.

We have that
\begin{align*}
 \|u-u_h\|_{L^2(\Omega)}^2 &= (u-u_h,\varphi) = (u-u_h, L_0^* \psi) \\
 &= a_{L_0}(t; u-u_h, \psi)  \\
 &= a_{L_0}(t; u-u_h, \psi-P_h^*(t,\psi)) + a_{L_0}(t; u-u_h, P_h^*(t,\psi))\\
&\leq  C_0 \|u-u_h\|_{H^1_0(\Omega)} \|\psi-P_h^*(t,\psi) \|_{H^1_0(\Omega)} \\
&\leq  C_2 h^{s-1} \|u\|_{H^s(\Omega)} h \|\psi \|_{H^2(\Omega)} \\
&\leq C_2  h^{s} \|u\|_{H^s(\Omega)} C_1 \|\varphi\|_{L^2(\Omega)}.
\end{align*}
Consequently,  we get that $\|u-P_h(t,u)\|_{L^2(\Omega)} \leq C h^{s} \|u\|_{H^s(\Omega)}$ for some constant $C>0$ independent of $t$ and $h$. The last result and inequality \eqref{eq:ineq_1} immediately lead to the desired result.

(ii)\quad By using the triangle inequality, it is sufficient to obtain appropriate bounds for $\|u - R_hu \|_{L^\infty(\Omega)}$ and  $\|R_hu-P_h(t,u)\|_{L^\infty(\Omega)}$. We recall that for polygonal domains, convex polyhedral domains, and smooth domains in $\mathbb{R}^N$  the Ritz projection $R_h$ is almost stable in $L^\infty(\Omega)$, that is
\[
\|R_h w\|_{L^\infty(\Omega)} \leq C_1\ell_h \|w\|_{L^\infty(\Omega)}\quad \forall w\in L^\infty(\Omega).
\]
We notice that $u-R_h u = (u-I_h u) -R_h(u-I_h u)$ where $I_h$ represents the interpolation operator from $C(\overline{\Omega})$ to $S_h^r(\Omega)$. Moreover, by using the properties of $I_h$ (see \cite{Livro2002_Ciarlet}) we get that
\begin{equation*}
\|u-R_h u\|_{L^\infty(\Omega)} \leq C_2\ell_h \|u-I_h u \|_{L^\infty(\Omega)} \leq C_3 \ell_h h^{s-\frac{n}{2}}\|u\|_{H^s(\Omega)},
\end{equation*}
where $C_2 = \max\{1,C_1\}$. On the other hand, using the  inverse inequality $\|\chi\|_{L^\infty(\Omega)} \leq C_4 h^{-\frac{n}{2}}\|\chi\|_{L^2(\Omega)}$ for $\chi\in S_h^r(\Omega)$ and property (i), we arrive at
\begin{align*}
\|R_hu-P_h(t,u)\|_{L^\infty(\Omega)} &\leq C_2 h^{-\frac{n}{2}} \|R_hu-P_h(t,u) \|_{L^2(\Omega)} \\
&\leq C_3 h^{s-\frac{n}{2}}\|u\|_{H^s(\Omega)}.
\end{align*} 

(iii)\quad Equation \eqref{eq:def_DP}  can be formally obtained from \eqref{eq:t-ritz_proj} applying the Leibniz rule. Moreover, as a direct consequence of the Lax-Milgram theorem and the properties of the form $a_{L_0}$ it follows by induction  that equation \eqref{eq:def_DP} defines functions $z_j :(0,T)\to S_h^r(\Omega))$. The proof that the functions $z_j$ coincide with the derivatives of $P_hu$ follows by induction after a lengthy but straightforward calculation, taking into account remark \ref{obs:prop_iii-iv}.

(iv)\quad We proceed by induction considering the base case where $j=0$; for this case the conclusion readily follows from property \eqref{eq:est_proj_i}. Assume that the estimate \eqref{eq:est_proj_iv} is valid for $j=0,\dots,q-1$, where $q-1<m$.

Let 
$w_h =  \partial_t^q (P_h u)(t) - P_h(t; \partial_t^q u(t))$, then from the definition of the Ritz projection and equation \eqref{eq:def_DP} with $j=q$,  we get that for all $\chi\in S_h^r(\Omega)$
\[
a_{L_0}(t; w_h, \chi) = \iim \sum_{l=0}^{q-1} \dbinom{p}{l} \gamma^{(j-l)}(t) (y\cdot\nabla 
[\partial_t^l u(t)-\partial_t^l (P_h u)(t)],\chi). 
\]
Setting $\chi = w_h$ and $e_h = u(t)-(P_hu)(t)$ we get that
\begin{align*}
 \|w_h\|_{H^1_0(\Omega)}^2 &\leq  2 Re\,\left\{ \iim \sum_{l=0}^{q-1} \dbinom{q}{l} \gamma^{(q-l)}(t) (y\cdot\nabla \partial_t^l e_h, w_h) \right\}\\
&\leq 2  Re\,\left\{\iim \sum_{l=0}^{q-1} \dbinom{q}{l} \gamma^{(q-l)}(t)  \left[ -n (\partial_t^l e_h, w_h) - (\partial_t^l e_h, y\cdot\nabla w_h) \right]\right\}\\
&\leq 2 \sum_{l=0}^{q-1} \dbinom{q}{l} \left|\gamma^{(q-l)}(t)\right|  \,\left| n (\partial_t^l e_h, w_h) + (\partial_t^l e_h, y\cdot\nabla w_h)\right|\\
&\leq 2 \sum_{l=0}^{q-1} \dbinom{q}{l} \left|\gamma^{(q-l)}(t)\right| \|\partial_t^l e_h\|_{L^2(\Omega)} \left(n \|w_h\|_{L^2(\Omega)} + d_\Omega \|w_h\|_{H_0^1(\Omega)}\right) \\
&\leq 2 (2^q-1) (n C_P + d_\Omega) \|\gamma\|_{C^q([0,T])}  \left(\max_{l=0,\dots,q-1} \|\partial_t^l e_h\|_{L^2(\Omega)} \right) \|w_h\|_{H^1_0(\Omega)},
\end{align*}
where $C_P$ represents the constant in the Poincaré inequality for functions in $H^1_0(\Omega)$. 
Consequently, from property (i) and the Poincaré inequality, it follows that for some constant $C>0$
\begin{equation}\label{eq:ineq_2}
 \max\left\{ \|w_h\|_{L^2(\Omega)}, \|w_h\|_{H^1_0(\Omega)} \right\}  \leq C h^{s}
\sum_{j=0}^{q-1}\|
\partial_t^j u(t)\|_{H^s(\Omega)}. 
\end{equation}
Finally, since 
$\partial_t^q u(t)-\partial_t^q (P_hu)(t) = \partial_t^q u(t) - P_h(t,\partial_t^q u(t)) - w_h$ using the triangle inequality and property \eqref{eq:est_proj_i}, we obtain the desired estimates.  
\end{proof}

\subsection{Error estimates for the semidiscrete problem}
%%%%%%%%%%%%%%%%%%%

This section aims to establish optimal error estimates for the solution to the semidiscrete problem \eqref{prob:aproximado}.  Our main result is given in the following theorem.

%========================================================================
\begin{thm}\label{teo:caso_semi-discreto}
%Suppose that $1\leq s\leq r$ if $d=1$ or $2\leq s\leq r$, if $d=2, 3$.   
Let $v$ and $v_h$ solve problems \eqref{prob:variational} and \eqref{prob:aproximado}, respectively. Assume that $v\in L^\infty(0,T; H^1_0(\Omega)\cap H^s(\Omega))$ and $\partial_{t} v \in L^2(0,T; H^1_0(\Omega)\cap H^s(\Omega))$ where $1\leq s\leq r$ for $n=1$, and $2\leq s\leq r$ for $n=2, 3$. If $\|P_hv(0)-v_{0h}\|_{L^2(\Omega)} \leq C_0 h^s$ for some constant $C_0>0$ then, there exist positive constants $h_0$, $C$ such that for $t\in [0,T]$ and $h\leq h_0$
\begin{equation*}
%\|v(t)-v_h(t)\|_{L^2(\Omega)} + h \|v(t)-v_h(t)\|_{H^1_0)(\Omega)} \leq Ch^s.
\|v(t)-v_h(t)\|_{L^2(\Omega)} + h \|v(t)-v_h(t)\|_{H^1_0(\Omega)}  \leq   Ch^s,
\end{equation*}
where the constant $C$ does not depend on $h$. 
\end{thm}
%========================================================================

Before starting the proof, we observe that the assumptions of the theorem imply that $v\in W^1(0,T; H^s(\Omega))$. 
\begin{proof}
The proof relies on a technique initiated by \cite{Artigo1973_Wheeler} for the analysis of parabolic problems. The main idea is to use the decomposition $v-v_h = \xi_h +\eta_h$  with $\xi_h = P_h v - v_h$ and $\eta_h = v-P_h v$,  get  appropriate bounds for $\xi_h$ and use properties of the Ritz projection to bound $\eta_h$. 

From equations \eqref{prob:variational} and \eqref{prob:aproximado},  we get that
\begin{equation}\label{eq:prob_dif}
(\partial_t \xi_h(t),\chi) = -\frac{\iim}{k^2(t)} a_L(t; \xi_h(t), \chi) - \frac{\iim}{k^2(t)} a_L(t; \eta_h(t), \chi)   %
- (g(v(t))-g(v_h(t)),\chi) - (\partial_t \eta_h(t),\chi)
\end{equation}
for any $\chi\in S_h^r(\Omega)$ and $t\leq T$ lying in the maximal  existence interval of the solution to the semidiscrete problem \eqref{prob:aproximado}. Next, by setting $\chi = \xi_h(t)$  in equation \eqref{eq:prob_dif} and noting that
\begin{equation*}
a_L(t; \eta_h(t),  \xi_h(t)) = a_{L_0}(t; \eta_h(t),  \xi_h(t)) - \lambda_0 (\eta_h(t),  \xi_h(t)) = - \lambda_0 (\eta_h(t),  \xi_h(t))
\end{equation*}
we get that
\begin{multline}\label{eq:est_xi1}
(\partial_t \xi_h(t),\xi_h(t)) = - \frac{\iim}{k^2(t)} a_L(t; \xi_h(t), \xi_h(t)) - \frac{\iim \lambda_0}{k^2(t)}(\eta_h(t), \xi_h(t)) \\- (g(v(t))-g(v_h(t)),\xi_h(t))- (\partial_t \eta_h(t),\xi_h(t)). 
\end{multline}

 We shall obtain an a priori estimate by taking the real part of \eqref{eq:est_xi1}. We have that
\begin{equation}
\begin{aligned}\label{eq:parcela1.1}
2 Re\,(\partial_t \xi_h(t), \xi_h(t)) &= (\partial_t \xi_h(t), \xi_h(t)) + \overline{(\partial_t \xi_h(t), \xi_h(t))} \\
&= (\partial_t \xi_h(t), \xi_h(t)) + ( \xi_h(t), \partial_t \xi_h(t)) \\
& = \frac{d}{dt}\|\xi_h(t)\|^2_{L^2(\Omega)}.
\end{aligned}
\end{equation}

Since
\[
\frac{\iim}{k^2(t)} a_L(t; \xi_h(t), \xi_h(t)) = \frac{\iim}{k^2(t)} \|\xi_h(t)\|^2_{H^1_0(\Omega)} - \frac{k'(t)}{k(t)} (y.\nabla\xi_h(t), \xi_h(t))
\]
and 
\[
\begin{aligned}
2Re\, (y\cdot\nabla\xi_h(t), \xi_h(t)) &= (y\cdot\nabla\xi_h(t), \xi_h(t)) +\overline{(y\cdot\nabla\xi_h(t), \xi_h(t))} \\
&= (y\cdot\nabla\xi_h(t), \xi_h(t)) +\overline{(y\cdot\nabla\xi_h(t), \xi_h(t))} \\
&= -(\xi_h(t), \nabla\cdot (y\xi_h(t))) +\overline{(y\cdot\nabla\xi_h(t), \xi_h(t))} \\
&= -(\xi_h(t), (\nabla\cdot y)\xi_h(t))) -(\xi_h(t), y\cdot \nabla \xi_h(t)) %\\&\qquad\qquad\qquad\qquad\qquad\quad 
+\overline{(y\cdot\nabla\xi_h(t), \xi_h(t))} \\
&= -n \|\xi_h(t)\|^2_{L^2(\Omega)} %
\end{aligned}
\]
we get that 
\begin{equation} \label{eq:parcela1.2}
2Re\, \frac{\iim}{k^2(t)} a_L(t; \xi_h(t), \xi_h(t)) =  \frac{n\,k'(t)}{k(t)} \|\xi_h(t)\|^2_{L^2(\Omega)}.
\end{equation}
Additionally, using  the Cauchy-Schwarz and Young inequalities, we get
\begin{equation}\label{eq:parcela1.3a}
\begin{aligned}
-2Re\, \frac{\iim \lambda_0}{k^2(t)}(\eta_h(t), \xi_h(t))  &\leq   \frac{2\lambda_0}{k^2(t)}\|\eta_h(t)\|_{L^2(\Omega)} \|\xi_h(t)\|_{L^2(\Omega)}\\
&\leq  \frac{\lambda_0}{k^2(t)}\left( \|\xi_h(t)\|^2_{L^2(\Omega)}  +\|\eta_h(t)\|^2_{L^2(\Omega)}  \right),
\end{aligned}
\end{equation}
\begin{equation}\label{eq:parcela1.3b}
\begin{aligned}
- 2Re\,(\partial_t \eta_h(t),\xi_h(t)) &\leq 2\|\partial_t \eta_h(t)\|_{L^2(\Omega)} \|\xi_h(t)\|_{L^2(\Omega)} \\
&\leq   \|\xi_h(t)\|^2_{L^2(\Omega)}  +\|\partial_t \eta_h(t)\|^2_{L^2(\Omega)}.  
\end{aligned}
\end{equation}

On the other hand, it is well-know (see, for instance, \cite{Cazenave2003}) that since $g(v) = |v|^\rho v$, there is a constant $c_\rho>0$ such that
\begin{equation}\label{eq:ineq_g}
|g(u_1)- g(u_2)|  \leq c_\rho(|u_1|^\rho + |u_2|^\rho)|u_1-u_2|.
\end{equation}
Thus, considering $u_1 = v(t)$ and $u_2=v_h(t)$ in the above equation, it would be helpful to have a bound for $|v(t)|^\rho + |v_h(t)|^\rho$ independent of $t$ and $h$. This will be achieved by employing a form of transfinite induction argument. %

Let 
\begin{equation}\label{eq:T_bar}
T^*_h = \sup\{ T^*\geq 0\,:\, \|v_h(t)\|_{L^\infty(\Omega)} \leq 2\|v\|_{L^\infty([0,T]\times\Omega)}\;\forall t\in[0,T^*]\}.
\end{equation} 
Since $v_{h0} = v_{h0}-P_hv(0) + P_hv(0) - v(0) + v(0)$, the triangle and inverse inequalities together with estimate \eqref{eq:est_proj_ii} and the hypotheses of the theorem yield
\begin{align*}
\|v_h(0)\|_{L^\infty(\Omega)} &\leq \|v_{h0}-P_hv(0)\|_{L^\infty(\Omega)} + \|P_hv(0) - v(0)\|_{L^\infty(\Omega)} 
+ \|v(0)\|_{L^\infty(\Omega)} \\
&\leq C_1 h^{-n/2} \|v_{h0}-P_hv(0)\|_{L^2(\Omega)} + C_2 \ell_h h^{s-n/2} \|v(0)\|_{H^s(\Omega)} % 
+ \|v\|_{L^\infty([0,T]\times\Omega)} \\
&\leq C_0 C_1 h^{s-n/2} + C_2 \ell_h h^{s-n/2} \|v(0)\|_{H^s(\Omega)}  + \|v\|_{L^\infty([0,T]\times\Omega)}.
\end{align*}
Moreover, since $s>n/2$ we can take $h'_0>0$ such that for $h\leq h'_0$ it holds that $\|v_h(0)\|_{L^\infty(\Omega)} < 2 \|v\|_{L^\infty([0,T]\times\Omega)}$. Therefore, the continuity of the semidiscrete solution $v_h$  implies that $T_h^* \in(0,T]$. From now on, we consider that $h\leq h'_0$.

Hence, for $t\in[0,T_h^*]$ using \eqref{eq:ineq_g} we obtain that 
\begin{equation}\label{eq:v_infty}
\begin{aligned} 
- 2Re\,(g(v(t))- g(v_h(t)),\xi_h(t))  &\leq 2C_\rho \left( \|\xi_h(t)+\eta_h(t)\|_{L^2(\Omega)} \right) \|\xi_h(t)\|_{L^2(\Omega)} \\
 &\leq C_\rho \left(3 \|\xi_h(t)\|^2_{L^2(\Omega)}+\|\eta_h(t)\|^2_{L^2(\Omega)} \right)
\end{aligned}
\end{equation}
where $C_\rho = c_\rho (1+2^\rho)\|v\|^\rho_{L^\infty([0,T]\times\Omega)}$
 
Using \eqref{eq:parcela1.1}--\eqref{eq:parcela1.3b} and \eqref{eq:v_infty}, it follows that
\begin{equation}
\frac{d}{dt}\|\xi_h(t)\|^2_{L^2(\Omega)}\leq C_1 \|\xi_h(t)\|^2_{L^2(\Omega)} + C_2 \|\eta_h(t)\|^2_{L^2(\Omega)} +    \|\partial_t \eta_h(t)\|^2_{L^2(\Omega)} 
\end{equation}
where 
\begin{align*}
C_1 &=  \max_{t\in[0,T]} | 1+3C_\rho-{n k'(t)}/{k(t)} + {\lambda_0}/{k^2(t)} | \\
C_2 &=   \max_{t\in[0,T]}  \left\{C_\rho +{\lambda_0}/{k^2(t)}\right\}.
\end{align*}
Thus, for $t\in[0,T_h^*]$  Lemma \ref{lem:ritz_proj} implies that 
\[
\frac{d}{dt}\|\xi_h(t)\|^2_{L^2(\Omega)}\leq C_1 \|\xi_h(t)\|^2_{L^2(\Omega)} + C_3 h^{2s} \left(\|v(t)\|^2_{H^s(\Omega)} + \|\partial_{t} v(t)\|^2_{H^s(\Omega)}\right)
\]
with a constant $C_3$ independent of $t$ and $h$. Consequently,
\[
\|\xi_h(t)\|^2_{L^2(\Omega)} \leq \|\xi_h(0)\|^2_{L^2(\Omega)} +C_3 h^{2s}\|v\|^2_{W^1(0,T;H^s(\Omega))} + C_1\int_0^t  \|\xi_h(t')\|^2_{L^2(\Omega)}\,dt',
\]
where 
\[
\|v\|^2_{W^1(0,T;H^s(\Omega))} = \int_0^T\left(\|v(t)\|^2_{H^s(\Omega)} + \|\partial_{t} v(t)\|^2_{H^s(\Omega)}\right) dt.
\]
Using the hypotheses of the theorem, we get
\[
\|\xi_h(t)\|^2_{L^2(\Omega)} \leq \left(C_0^2 +C_3 \|v\|^2_{W^1(0,T;H^s(\Omega))}\right)h^{2s} + C_1\int_0^t  \|\xi_h(t')\|^2_{L^2(\Omega)}\,dt'
\]
Hence, by applying the Gronwall Lemma, we get for $t\in[0,T_h^*]$
\begin{equation}\label{eq:est_xi}
\|\xi_h(t)\|_{L^2(\Omega)} \leq  C_4 h^{s}
\end{equation}
where $C_4^2 = e^{C_1T}\,(C_0^2 +C_3 \|v\|^2_{W^1(0,T;H^s(\Omega))} )$. Moreover, the inverse inequality implies that 
\[
\|\xi_h(t)\|_{H^1_0(\Omega)} \leq  C_5 h^{s-1}.
\]
Finally, using the triangle inequality and property \eqref{eq:est_proj_i}, it follows that for $t\in[0,T_h^*]$
\begin{equation}\label{eq:est_error}
\|v(t) - v_h(t)\|_{L^2(\Omega)} + h\|v(t) - v_h(t)\|_{H^1_0(\Omega)}\leq  C h^{s}
\end{equation}
where the constant $C$ is independent of $t$ and $h$.

In order to finish the proof of the theorem, we need to remove the restriction $t\in[0,T_h^*]$. Observe that $v_h = v -\xi_h -\eta_h$. Then, for $t\in[0,T_h^*]$ using the triangle and inverse inequalities and property \eqref{eq:est_proj_ii},  we get
\begin{align*}
\|v_h(t)\|_{L^\infty(\Omega)} &\leq \|v(t)\|_{L^\infty(\Omega)} + \|\xi_h(t)\|_{L^\infty(\Omega)} + \|\eta_h(t)\|_{L^\infty(\Omega)} \\
&\leq \|v(t)\|_{L^\infty(\Omega)} + C_6 h^{-n/2} \|\xi_h(t)\|_{L^2(\Omega)} 
%\\&\qquad\qquad\qquad\qquad\qquad 
+ C_7\ell_h h^{s-n/2} \|v\|_{L^\infty(0,T; H^s(\Omega))},
\end{align*}
with positive constants $C_6$ and $C_7$ independent of $h$. 
%We set 
%\[
%M_h(t) = C_6 h^{-n/2} \|\xi_h(t)\|_{L^2(\Omega)} + C_7\ell_h h^{s-n/2} \|v\|_{L^\infty(0,T; H^s(\Omega))}
%\]
%and apply estimate \eqref{eq:est_xi} to obtain
%\[ 
%M_h(t) \leq  h^{s-n/2} (C_8  + C_7\ell_h  \|v\|_{L^\infty(0,T; H^s(\Omega))} ).
%\]
%%\end{align*}
%Since $s-n/2>0$, we can take $h_0$ such that if  $h<h_0$ then $M_h(t) < \|v\|_{L^\infty([0,T]\times \Omega)}$ for $t\in[0,T_h^*]$. 
We apply estimate \eqref{eq:est_xi} to obtain
\[ 
\|v_h(t)\|_{L^\infty(\Omega)} \leq \|v\|_{L^\infty([0,T]\times \Omega)} + h^{s-n/2} (C_8  + C_7\ell_h  \|v\|_{L^\infty(0,T; H^s(\Omega))} ).
\]
%\end{align*}
Since $s-n/2>0$, we can take $h_0 \leq h'_0$ such that if  $h<h_0$ then $\|v_h(t)\|_{L^\infty(\Omega)} < 2\|v\|_{L^\infty([0,T]\times \Omega)}$ for $t\in[0,T_h^*]$. 

Consequently, if $T_h^* < T$, then the approximate solution $v_h$ would exist in some interval $I$ larger than $[0, T_h^*]$.  Therefore, since $v_h$ is continuous, there is a $\delta>0$ such that  
$\|v_h(t)\|_{L^\infty(\Omega)}< 2\|v\|_{L^\infty([0,T]\times \Omega)}$  for $t\in[T_h^*, T_h^*+\delta)$. This is a contradiction with the definition of $T_h^*$, and we conclude that $T_h^*=T$ whenever $h<h_0$; hence estimate \eqref{eq:est_error} holds for all $t\in[0,T]$. 
\end{proof}

\subsection{Error estimates for the fully discrete problem}
%%%%%%%%%%%%%%%%%%%

For the error estimate of the solution to the fully discrete problem, we assume that the solution to \eqref{prob:variational} is more regular than that of theorem \ref{teo:caso_semi-discreto} for the semidiscrete case. We shall consider that the solution $v$ of problem \eqref{prob:variational} satisfies the following assumptions.
\begin{enumerate}
\item[\textbf{A1:}] $v\in L^\infty(0,T; H^s(\Omega) \cap H^1_0(\Omega))$ and $\partial_{t} v \in L^2(0,T; H^s(\Omega) \cap H^1_0(\Omega))$
\item[\textbf{A2:}] $v, \partial_{t} v, \partial_{t}^2 v \in L^2(0,T; H^2(\Omega))$, i.e. $v\in W^2(0,T; H^2(\Omega))$
\item[\textbf{A3:}] $\partial_t^3 v \in L^2(0,T; L^2(\Omega))$
\item[\textbf{A4:}] $k\in C^4([0,T])$, $k(t)\geq k_0>0$ for $t\in[0,T]$
\end{enumerate}

Our main result is given in the following theorem.

%========================================================================
\begin{thm}\label{teo:caso_discreto}
%Suppose that $1\leq s\leq r$ if $d=1$ or $2\leq s\leq r$, if $d=2, 3$.   
Let $v$ and $\{U^m\}$ (m=1,2,\dots,N) solve problems \eqref{prob:variational} and \eqref{prob_aprox:discrete_all}, respectively. Assume that $v$ satisfies assumptions \textbf{A1--A4}  where $1\leq s\leq r$ for $n=1$, and $2\leq s\leq r$ for $n=2, 3$.  %Suppose that $\|P_hv(0)-U^0\|_{L^2(\Omega)} \leq C_0 h^s$ and $\tau \leq C_1 h^{{s}/{2}}$ for some constants $C_0,\,C_1>0$ then, there exist positive constants $h_0$ and $C$ such that for  $h\leq h_0$
Suppose that  $\|P_hv(0)-U^0\|_{L^2(\Omega)} \leq C_0 h^s$ for some constant $C_0>0$ and $\tau = {o}(h^{n/4})$, then there exist positive constants $h_0$ and $\tau_0$ such that if $h\leq h_0$ and $\tau \leq  \tau_0$ then
\begin{equation*}
%\|v(t)-v_h(t)\|_{L^2(\Omega)} + h \|v(t)-v_h(t)\|_{H^1_0)(\Omega)} \leq Ch^s.
\max_{m= 0,\dots, T/\tau} \|v(t_m)-U^m\|_{L^2(\Omega)}   \leq   C(\tau^2 + h^s),
\end{equation*}
where the constant $C$ does not depend on $h$ or $\tau$. 
%% Remark. we have that $v\in C([0,T];H^s(\Omega))$
\end{thm}
%========================================================================

\begin{proof}
We proceed as in theorem \ref{teo:caso_semi-discreto} by taking advantage of the estimates associated with the Ritz projection presented in lemma \ref{lem:ritz_proj}.

We apply the decomposition $v(t_m)-U^m = \theta^m +\eta^m$ where $\theta^m = (P_h v)(t_m) - U^m$ and $\eta^m = v(t_m)-(P_h v)(t_m)$. Using similar definitions to that in \eqref{eq:discrete_diff} and \eqref{eq:notation_all}, we also consider the quantities 
$\delta_\tau \theta^m$, $\delta_\tau \eta^m$, $\hat{\theta}^m$, $\hat{\eta}^m$, $\tilde{\theta}^m$, $\tilde{\eta}^m$.  Our goal is to obtain appropriate estimates for the projected errors $\theta^m$. This will be achieved by establishing a recursive relation and applying a discrete Gronwall lemma.   

From \eqref{prob:variational} and \eqref{prob_aprox:discrete_all}, for any $\chi\in S^r_h(\Omega)$ we get 
\begin{align}\label{eq:approx_discrete}
& (\delta_\tau \theta^m, \chi)
= (\delta_\tau v^m - \partial_{t} v(t_{m-\frac{1}{2}}), \chi)
- (\delta_\tau \eta^m, \chi)
+ \frac{\iim}{k^2_{m-\frac{1}{2}}} a_{m-\frac{1}{2}}(\hat{v}^m - v(t_{m-\frac{1}{2}}), \chi)
\nonumber\\ 
& \quad
+ \frac{\iim}{k^2_{m-\frac{1}{2}}} a_{m-\frac{1}{2}}( \eta^{m-\frac{1}{2}} - \hat{\eta}^m, \chi)
- \frac{\iim}{k^2_{m-\frac{1}{2}}} a_{m-\frac{1}{2}}( \eta^{m-\frac{1}{2}}, \chi) 
\nonumber\\ 
& \quad
- \frac{\iim}{k^2_{m-\frac{1}{2}}} a_{m-\frac{1}{2}}( \hat{\theta}^m, \chi)   
- (g(v(t_{m-\frac{1}{2}})) - g^m ,\chi),
\end{align}
where $a_{m-\frac{1}{2}}( \cdot, \cdot)= a_L(t_{m-\frac{1}{2}}; \cdot, \cdot)$. We set  $\chi=2\hat{\theta}^m$ and take the real part of \eqref{eq:approx_discrete}. For the left-hand side term, we get that 
\begin{equation}
Re\, (\delta_\tau \theta^m, 2\hat{\theta}^m) = \frac{1}{\tau}\left( \|\theta^m\|^2_{L^2(\Omega)} - \|\theta^{m-1}\|^2_{L^2(\Omega)} \right).
\end{equation}

Our next goal is to bound the real part of the seven terms appearing on the right-hand side of \eqref{eq:approx_discrete}. 

For $\omega^m_1 =\delta_\tau v^m - \partial_{t} v(t_{m-\frac{1}{2}})$ we have 
\[
\begin{aligned}
\omega^m_1 &= \frac{1}{\tau} \int_{t_{m-1}}^{t_{m}} (\partial_{t} v(t) - \partial_{t} v(t_{m-\frac{1}{2}}) dt \\
&= \frac{1}{\tau}\int_{t_{m-1}}^{t_{m}}\left( (t-t_{m-\frac{1}{2}})\partial_t^2 v(t_{m-\frac{1}{2}}) + \int_{t_{m-\frac{1}{2}}}^t (t-l) \partial_t^3 v(l) dl \right)dt \\
&= \frac{1}{\tau}\int_{t_{m-1}}^{t_{m}} \int_{t_{m-\frac{1}{2}}}^t (t-l) \partial_t^3 v(l) dl\,dt \\
&= \frac{1}{2\tau}\left( \int_{t_{m-1}}^{t_{m-\frac{1}{2}}} (t-t_{m-1})^2 \partial_t^3 v(t)\,dt +  \int_{t_{m-\frac{1}{2}}}^{t_{m}} (t_{m}-t)^2 \partial_t^3 v(t)\,dt \right)
\end{aligned}
\]
 then 
\begin{equation*}%\label{eq:omega_1}
\begin{aligned}
\|\omega^m_1\|_{L^2(\Omega)} &\leq\frac{1}{2\tau}\left( \frac{\tau^2}{4} \int_{t_{m-1}}^{t_{m-\frac{1}{2}}} \|  \partial_t^3 v(t)\|_{L^2(\Omega)}\,dt +  \frac{\tau^2}{4} \int_{t_{m-\frac{1}{2}}}^{t_{m}} \|  \partial_t^3 v(t)\|_{L^2(\Omega)}\,dt  \right) \\
&\leq \frac{\tau}{8} \int_{t_{m-1}}^{t_{m}} \|  \partial_t^3 v(t)\|_{L^2(\Omega)}\,dt 
\end{aligned}
\end{equation*}
and
\begin{equation}\label{eq:S1}
\begin{aligned}
S^m_1 &= Re\, (\delta_\tau v^m - \partial_{t} v(t_{m-\frac{1}{2}}), 2\hat{\theta}^m) \\
&= Re\, (\omega^m_1, 2\hat{\theta}^m) \\
& \leq  \|\omega^m_1\|_{L^2(\Omega)} \| 2\hat{\theta}^m\|_{L^2(\Omega)} \\
&\leq \frac{\tau}{8} \left( \int_{t_{m-1}}^{t_{m}} \|  \partial_t^3 v(t)\|_{L^2(\Omega)}\,dt \right)\left(\| \theta^m\|_{L^2(\Omega)} +  \| \theta^{m-1}\|_{L^2(\Omega)} \right).
\end{aligned}
\end{equation}

Using property \eqref{eq:est_proj_iv} we have 
\[
\begin{aligned}
\|\delta_\tau \eta^m\|_{L^2(\Omega)} &= \frac{1}{\tau}\left\|\int_{t_{m-1}}^{t_m} \partial_t \eta(t) dt \right\|_{L^2(\Omega)} \leq \frac{1}{\tau}\int_{t_{m-1}}^{t_m}  \|\partial_t \eta(t) \|_{L^2(\Omega)} dt\\
&\leq \frac{C_1h^{s}}{\tau} \int_{t_{m-1}}^{t_m}  \left( \|v(t) \|_{H^s(\Omega)}  +\|\partial_{t} v(t) \|_{H^s(\Omega)} \right) dt  %
\end{aligned}
\]
then 
\begin{equation}\label{eq:S2}
\begin{aligned}
S^m_2 &= -Re\, (\delta_\tau \eta^m, 2\hat{\theta}^m) \leq \|\delta_\tau \eta^m\|_{L^2(\Omega)} \|2\hat{\theta}^m\|_{L^2(\Omega)} \\
&\leq  \frac{C_1 h^{s}}{\tau} \left[\int_{t_{m-1}}^{t_m}  \left( \|v(t) \|_{H^s(\Omega)}  +\|\partial_{t} v(t) \|_{H^s(\Omega)} \right)dt \right]  %\\&\qquad \qquad \qquad \qquad 
\times \left( \| \theta^m\|_{L^2(\Omega)} + \| \theta^{m-1}\|_{L^2(\Omega)} \right).
\end{aligned}
\end{equation}

 Set 
 \begin{equation}\label{eq:omega2^m}
\omega^m_2 =\hat{v}^m - v(t_{m-\frac{1}{2}}),
\end{equation}
 then we have that
\begin{equation}\label{eq:S3_a}
 \begin{aligned}
 S^m_3 &= Re\,  \frac{\iim}{k^2_{m-\frac{1}{2}}} a_{m-\frac{1}{2}}( \hat{v}^m - v(t_{m-\frac{1}{2}}), 2\hat{\theta}^m) \\
 &\leq \frac{1}{k^2_{m-\frac{1}{2}}} |a_{m-\frac{1}{2}}( \omega^m_2, 2\hat{\theta}^m)| \\
 &\leq \frac{1}{k^2_{m-\frac{1}{2}}} |-(\Delta \omega^m_2, 2\hat{\theta}^m) +\iim k'_{m-\frac{1}{2}} k_{m-\frac{1}{2}}( y\cdot\nabla\omega^m_2, 2\hat{\theta}^m)| \\
 & \leq \frac{1}{k^2_{m-\frac{1}{2}}} \left( \|\Delta \omega^m_2\|_{L^2(\Omega)} +  |k'_{m-\frac{1}{2}} k_{m-\frac{1}{2}}|\,d_\Omega  \|\omega^m_2\|_{H^1_0(\Omega)} \right) \|2\hat{\theta}^m\|_{L^2(\Omega)}.
\end{aligned}
\end{equation}
Since 
\[
\begin{aligned}
\omega^m_2 &= \frac{1}{2} \left( v(t_m) - v(t_{m-\frac{1}{2}}) + v(t_{m-1}) - v(t_{m-\frac{1}{2}})  \right) \\
&= \frac{1}{2} \left( \int_{t_{m-\frac{1}{2}}}^{t_{m}} (t_m-t) \partial_t^2 v(t)  dt + \int_{t_{m-\frac{1}{2}}}^{t_{m-1}} (t_{m-1}-t) \partial_t^2 v(t)  dt   \right) 
\end{aligned}
\]
we get 
\[
\begin{aligned}
\|\Delta\omega^m_2\|_{L^2(\Omega)} &\leq \frac{1}{2} \left( \frac{\tau}{2}\int_{t_{m-\frac{1}{2}}}^{t_{m}} \|\Delta \partial_t^2 v(t)\|_{L^2(\Omega)}  dt + \frac{\tau}{2} \int_{t_{m-1}}^{t_{m-\frac{1}{2}}} \| \Delta \partial_t^2 v(t)\|_{L^2(\Omega)}  dt   \right) \\
&\leq \frac{\tau}{4}  \int_{t_{m-1}}^{t_{m}} \| \Delta \partial_t^2 v(t)\|_{L^2(\Omega)}  dt   \\
&\leq \frac{\tau}{4}   \int_{t_{m-1}}^{t_{m}} \|\partial_t^2 v(t)\|_{H^2(\Omega)}  dt 
\end{aligned}
\]
and 
\[
\|\omega^m_2\|_{H^1_0(\Omega)} \leq \frac{\tau}{4}  \int_{t_{m-1}}^{t_{m}} \|\partial_t^2 v(t)\|_{H^1_0(\Omega)}  dt.
\]
Consequently,
\begin{equation}\label{eq:S3}
S^m_3 \leq  C_2 \tau \left( \int_{t_{m-1}}^{t_{m}} \|\partial_t^2 v(t)\|_{H^2(\Omega)}  dt \right)
\left( \| \theta^m\|_{L^2(\Omega)} +  \| \theta^{m-1}\|_{L^2(\Omega)} \right)
\end{equation}
where $C_2 = \max_{t\in[0,T]}\left\{ \frac{1}{2k^2(t)},  \frac{d_\Omega|k'(t)|}{2k(t)} \right\}$.

Proceeding, as in the previous case and applying property \eqref{eq:est_proj_iv}, we get for $\omega^m_3=\eta^{m-\frac{1}{2}} - \hat{\eta}^m$
\begin{align*}
\|\omega^m_3\|_{H^1_0(\Omega)} &\leq \frac{\tau}{4}   \left( \int_{t_{m-1}}^{t_{m}} \|\partial_t^2 \eta(t)\|_{H^1_0(\Omega)}  dt \right)\\
&\leq C_3 h \tau    \int_{t_{m-1}}^{t_{m}}  \left( \|v\|_{H^2(\Omega)}+  \|\partial_{t} v\|_{H^2(\Omega)}+ \|\partial_t^2 v(t)\|_{H^2(\Omega)} \right) dt .
\end{align*}
Moreover, using the inverse inequality we get
\begin{equation}\label{eq:S4}
\begin{aligned}
 S^m_4 &= Re\,\frac{\iim}{k^2_{m-\frac{1}{2}}} a_{m-\frac{1}{2}}( \omega^m_3, 2\hat{\theta}^m) \leq \frac{1}{k^2_{m-\frac{1}{2}}} \|\omega^m_3\|_{H^1_0(\Omega)} \| 2\hat{\theta}^m\|_{H^1_0(\Omega)} \\
 & \leq \frac{1}{ h k^2_{m-\frac{1}{2}}} \|\omega^m_3\|_{H^1_0(\Omega)} \| 2\hat{\theta}^m\|_{L^2(\Omega)} \\ 
 &\leq C_4 \tau  \left[\int_{t_{m-1}}^{t_{m}}  \left( \|v\|_{H^2(\Omega)}+  \|\partial_{t} v\|_{H^2(\Omega)}+ \|\partial_t^2 v(t)\|_{H^2(\Omega)} \right)dt\right] %\\ &\qquad\qquad\qquad \qquad\qquad 
 \times \left( \| \theta^m\|_{L^2(\Omega)} +  \| \theta^{m-1}\|_{L^2(\Omega)} \right)
\end{aligned}
\end{equation}
where $C_4 = C_3\max_{t\in[0,T]} \frac{1}{k^2(t)}.$

We also have that 
\[
S^m_5 = -  Re\,\frac{\iim}{k^2_{m-\frac{1}{2}}} a_{m-\frac{1}{2}}( \eta^{m-\frac{1}{2}}, 2\hat{\theta}^m) =  Re\,\frac{i\lambda_0}{k^2_{m-\frac{1}{2}}} ( \eta^{m-\frac{1}{2}}, 2\hat{\theta}^m),
\]
and property \eqref{eq:est_proj_i} gives that
\begin{align*}
\| \eta^{m-\frac{1}{2}}\|_{L^2(\Omega)} &\leq C_5 h^s \|v(t_{m-\frac{1}{2}})\|_{H^s(\Omega)}.
%\\&\leq C_6 h^s \|v\|_{W^1(0,T; H^1_0(\Omega)\cap H^s(\Omega))} .
\end{align*}
We conclude that
\begin{equation}\label{eq:S5}
S^m_5 \leq  C_6 h^{s}  \|v\|_{L^\infty(0,T;H^s(\Omega))} \left( \| \theta^m\|_{L^2(\Omega)} +  \| \theta^{m-1}\|_{L^2(\Omega)} \right) ,
\end{equation}
where $C_6 = C_5 \lambda_0 \max_{t\in[0,T]} \frac{1}{k^2(t)}$. %\max_{t\in[0,T]} \frac{1}{k^2(t)}$.% 

Applying the same reasoning  as in \eqref{eq:parcela1.2} we get 
\begin{equation}\label{eq:S6}
\begin{aligned}
S^m_6 &= - Re\, \frac{\iim}{k^2_{m-\frac{1}{2}}} a_{m-\frac{1}{2}}( \hat{\theta}^m, 2\hat{\theta}^m) = - \frac{n k'_{m-\frac{1}{2}}}{k_{m-\frac{1}{2}}} \|\hat{\theta}^m\|^2_{L^2(\Omega)} \\
&\leq  C_8 \left( \| \theta^m\|^2_{L^2(\Omega)} +  \| \theta^{m-1}\|^2_{L^2(\Omega)}  \right)^2
\end{aligned}
\end{equation}
where $C_8 = n \max_{t\in[0,T]} \frac{ |k'(t)|}{4k(t)}$.

%Now, to bound the term $S^m_7 = -Re\,(g(v(t_{m-\frac{1}{2}})) - g^m , 2\hat{\theta}^m)$ we proceed as in the proof of theorem \ref{teo:caso_semi-discreto}.   
In order to bound the terms $S^m_7 = -Re\,(g(v(t_{m-\frac{1}{2}})) - g^m , 2\hat{\theta}^m)$, we need to consider three different cases for the index $m$ (see \eqref{eq:notation_all}). %and proceed as in the proof of theorem \ref{teo:caso_semi-discreto}. 
%%%

We start with the case $m=1^-$. We have that $g^{1^-} = g(U^0)$, and from inequality \eqref{eq:ineq_g}  one gets that 
\begin{equation}\label{eq:ineq_g_c1}
|g(v(t_{\frac{1}{2}})) - g^{1^-})|  \leq c_\rho(|v(t_{\frac{1}{2}})|^\rho + |{U}^0|^\rho) |v(t_{\frac{1}{2}}) - {U}^0|.
\end{equation}
Next, by taking $h$ sufficiently small, it follows that  
\begin{equation}\label{eq:U0inf_b}
\|U^{0}\|_{L^\infty(\Omega)} \leq 2\|v\|_{L^\infty([0,T]\times\Omega)}.
\end{equation}
 Indeed, from property \eqref{eq:est_proj_ii}, the triangle and inverse inequalities and the hypotheses of the theorem, we get
\[
\begin{aligned}
\|U^0\|_{L^\infty(\Omega)} &\leq \|U^{0} - P_hv(0)\|_{L^\infty(\Omega)} + \|P_hv(0) - v(0)\|_{L^\infty(\Omega)} +\|v(0)\|_{L^\infty(\Omega)} \\
&\leq C_ 0 h^{s-n/2} + C_9 \ell_h h^{s-n/2} \|v(0)\|_{H^s(\Omega)} +\|v(0)\|_{L^\infty(\Omega)}.
\end{aligned}
\]
Using that $s-n/2>0$ and choosing $h<h'_0$ with a sufficiently small $h'_0$, we arrive at the desired conclusion.

Applying \eqref{eq:U0inf_b} into \eqref{eq:ineq_g_c1}, we get that  
\begin{equation}\label{eq:ineq_g2}
|(g(v(t_{\frac{1}{2}})) - g({U}^0) ,\hat{\theta}^m)| \leq C_{10} \| v(t_{\frac{1}{2}}) - {U}^0 \|_{L^2(\Omega)} \|\hat{\theta}^{1^-}\|_{L^2(\Omega)}
\end{equation}
where $C_{10} = c_\rho (1+4^\rho)\|v\|_{L^\infty([0,T]\times\Omega)}$.
Consider the decomposition  
\[
v(t_{\frac{1}{2}}) - {U}^0 = \omega_4 + {\eta}^0 + {\theta}^0
\] 
where $\omega_4 = v(t_{\frac{1}{2}}) - v(0)$.
We have that
\begin{equation}\label{eq:ineq_omega4}
\|\omega_4 \|_{L^2(\Omega)}  = \left\|\int_{0}^{\frac{\tau}{2}} \partial_{t} v(t)  dt \right\|_{L^2(\Omega)} \leq  \frac{\tau}{2}  \| \partial_{t} v\|_{L^\infty(0,T; L^2(\Omega))}, 
\end{equation}
and from property \eqref{eq:est_proj_i} we get
\begin{equation}\label{eq:ineq_eta0}
\|{\eta}^{0}\|_{L^2(\Omega)}  \leq C_{11} h^s  \|v(0)\|_{H^s(\Omega)} 
 \leq C_{11} h^s \|v\|_{L^\infty(0,T; H^s(\Omega))}.
\end{equation}
Hence, we have
 \begin{equation}\label{eq:S7_1-}
 S^{1^-}_7  \leq 
 C_{12}  \Bigl( \tau  \| \partial_{t} v \|^2_{L^\infty(0,T; L^2(\Omega))} + h^s \|v\|_{L^\infty(0,T; H^s(\Omega))} 
 +  \|\theta^0\|_{L^2(\Omega)}\Bigr) \\
 \times\left( \| \theta^{1^-}\|_{L^2(\Omega)} +  \| \theta^{0}\|_{L^2(\Omega)}  \right).
 \end{equation}
 
 Collecting all the bounds from equations \eqref{eq:S1}, \eqref{eq:S2}, \eqref{eq:S3}, \eqref{eq:S4}, \eqref{eq:S5}, \eqref{eq:S6} and \eqref{eq:S7_1-}, we arrive at the conclusion that 
\begin{equation}\label{eq:ineq_1-a}
\|\theta^{1^-}\|_{L^2(\Omega)} \leq  A_{1}^- \tau^2 + B_1^- h^{s} + \left(1+C_{1}^-\tau \right) \|\theta^{0}\|_{L^2(\Omega)}  + \bar{D} \tau   \| \theta^{1^-}\|_{L^2(\Omega)} 
\end{equation}
where
\begin{equation*}
\begin{split}
&A_1^- =  C_{13} \Bigl(\int_{t_{0}}^{t_{1}} \bigl (\|  v(t)\|_{H^2(\Omega)}+ \|  \partial_{t} v(t)\|_{H^2(\Omega)}+ \|  \partial_t^2 v(t)\|_{H^2(\Omega)}  
+  \| \partial_t^3 v(t)\|_{L^2(\Omega)} \bigr) dt  %
+ \|\partial_{t} v\|_{L^\infty(0,T; L^2(\Omega))}\Bigr), \\
& B_1^- = C_{14}\Bigl( \tau \|v\|_{L^\infty(0,T;H^s(\Omega))}  
+ \int_{t_0}^{t_1}  \left( \|v(t) \|_{H^s(\Omega)}  +\|\partial_{t} v(t) \|_{H^s(\Omega)} \right) dt \Bigr), \\
& C_1^- = C_{12} +  \bar{D}, \\
& \bar{D} = n\max_{t\in[0,T]} \frac{ |k'(t)|}{4k(t)}. %%
\end{split}
\end{equation*}
Moreover, let $\tau_0=\frac{1}{2\bar{D}}$ then if $\tau \leq \tau_0$ we have that %
\begin{equation}%\label{eq:ineq_1-b}
\|\theta^{1^-}\|_{L^2(\Omega)} \leq  2\left(A_{1}^- \tau^2 + B_1^- h^{s}  + (1+C_1^-\tau) \|\theta^{0}\|_{L^2(\Omega)}  \right),
\end{equation}
and consequently
\begin{equation}\label{eq:ineq_1-b}
\|\theta^{1^-}\|_{L^2(\Omega)} \leq  E_{1}^- ( \tau^2 + h^{s} )
\end{equation}
where $E_{1}^- = \max\{2A_{1}^-, 3C_0 + 2B_1^- \}$.

Next, we consider the case $m=1$ and will proceed as in the previous case. 
By taking $h$ and $\tau$ sufficiently small, we have % and $\tau\leq\min\{\tau_0,h^q\}$, we have  
\begin{equation}\label{eq:U1-inf_b}
\|U^{1^-}\|_{L^\infty(\Omega)} \leq 2\|v\|_{L^\infty([0,T]\times\Omega)}.
\end{equation}
Indeed, using property \eqref{eq:est_proj_ii}, the triangle and inverse inequalities,   estimate \eqref{eq:ineq_1-b} and the hypotheses of the theorem, we have that
\[
\begin{aligned}
\|U^{1^-}\|_{L^\infty(\Omega)} &\leq \|\theta^{1^-} \|_{L^\infty(\Omega)} + \|\eta^1\|_{L^\infty(\Omega)} +\|v(t_1)\|_{L^\infty(\Omega)} \\
&\leq C_ {14} h^{-n/2}\|\theta^{1^-} \|_{L^2(\Omega)} + C_{15} \ell_h h^{s-n/2} \|v(t_1)\|_{H^s(\Omega)}  +\|v(t_1)\|_{L^\infty(\Omega)}\\
&\leq 2 C_ {14}E_{1}^-  (\tau^2 h^{-n/2} + h^{s-n/2}) %\\&\qquad \qquad\qquad 
+ C_{15} \ell_h h^{s-n/2} \|v\|_{L^\infty(0,T;H^s(\Omega))}  +\|v\|_{L^\infty([0,T]\times\Omega)}. %
\end{aligned}
\]
Noting that $s-n/2>0$ and $\tau^2 = o(h^{n/2})$ we can choose $h''_0>0$ and $\tau_0'>0$ such that for any $h\leq h''_0$
\[
2 C_ {14}E_{1}^-  (\tau^2 h^{-n/2} + h^{s-n/2}) + C_{15} \ell_h h^{s-n/2} \|v\|_{L^\infty(0,T;H^s(\Omega))} 
\leq \|v\|_{L^\infty([0,T]\times\Omega)}
\]
which proves the required inequality. 

Since $g^1 = g(\hat{U}^{1^-})$ applying inequalities \eqref{eq:ineq_g} and \eqref{eq:U1-inf_b},  we get that  
\begin{equation}\label{eq:ineq_gc12}
|(g(v(t_{\frac{1}{2}})) - g^1 ,\hat{\theta}^1)| \leq C_{16} \| v(t_{\frac{1}{2}}) - \hat{U}^{1^-} \|_{L^2(\Omega)} \|\hat{\theta}^{1}\|_{L^2(\Omega)}
\end{equation}
where $C_{16} = c_\rho (1+4^\rho)\|v\|_{L^\infty([0,T]\times\Omega)}$.
Consider the decomposition  
\[
v(t_{\frac{1}{2}}) - \hat{U}^{1^-} = -\omega_2^1 + \hat{\eta}^1 + \hat{\theta}^{1^-}
\] 
where $\omega_2^1 =  \hat{v}^1 - v(t_{\frac{1}{2}})$ as defined in \eqref{eq:omega2^m}.
We have that
\begin{equation}\label{eq:ineq_omega2^1}
\|\omega^1_2\|_{L^2(\Omega)} \leq \frac{\tau}{4}  \int_{t_{0}}^{t_{1}} \|\partial_t^2 v(t)\|_{L^2(\Omega)}  dt, 
\end{equation}
and from property \eqref{eq:est_proj_i} we get that
\begin{equation}\label{eq:ineq_hat_eta1}
\begin{aligned}
\|\hat{\eta}^{1}\|_{L^2(\Omega)}  &\leq C_{17} h^s  \left(\frac{\|v(t_0)\|_{H^s(\Omega)} + \|v(t_1)\|_{H^s(\Omega)}}{2} \right)\\
 & \leq C_{17} h^s \|v\|_{L^\infty(0,T; H^s(\Omega))}.
\end{aligned}
\end{equation}
Hence, we have
\begin{multline}\label{eq:S7_1}
 S^{1}_7  \leq 
 C_{18}  \Bigl(  \tau  \int_{t_0}^{t_1} \| \partial_t^2 v(t) \|_{L^2(\Omega)} dt + h^s \|v\|_{L^\infty(0,T; H^s(\Omega))}  
 +  \|\theta^0\|_{L^2(\Omega)} +  \|\theta^{1^-}\|_{L^2(\Omega)} \Bigr)
 \\ \times \left( \| \theta^{1}\|_{L^2(\Omega)} +  \| \theta^{0}\|_{L^2(\Omega)}  \right).
 \end{multline}
 
 Finally, using  estimates \eqref{eq:S1}, \eqref{eq:S2}, \eqref{eq:S3}, \eqref{eq:S4}, \eqref{eq:S5}, \eqref{eq:S6} and \eqref{eq:S7_1}, we get that 
\begin{equation}\label{eq:ineq_1a}
\|\theta^{1}\|_{L^2(\Omega)} \leq \bar{A}_{1} \tau^2 + \bar{B}_1 h^{s}+ (1+\bar{C}_1\tau) \|\theta^{0}\|_{L^2(\Omega)} 
+ C_{18}\tau \|\theta^{1^-}\|_{L^2(\Omega)} + \bar{D} \tau  \| \theta^{1}\|_{L^2(\Omega)} 
\end{equation}
where
\begin{equation}
\begin{split}
\bar{A}_1 &=  C_{19} \int_{t_{0}}^{t_{1}} \bigl (\|  v(t)\|_{H^2(\Omega)}+ \|  \partial_{t} v(t)\|_{H^2(\Omega)}+ \|  \partial_t^2 v(t)\|_{H^2(\Omega)} 
 +  \| \partial_t^3 v(t)\|_{L^2(\Omega)} \bigr) dt,  \\
\bar{B}_1 &= C_{20}\Bigl( \tau \|v\|_{L^\infty(0,T;H^s(\Omega))}  + \int_{t_0}^{t_1}  \left( \|v(t) \|_{H^s(\Omega)}  +\|\partial_{t} v(t) \|_{H^s(\Omega)} \right) dt \Bigr), \\
\bar{C}_1 &= C_{18} +  \bar{D}.
\end{split}
\end{equation}

Moreover, for $\tau \leq \tau_0$ we have that %
\begin{equation}%\label{eq:ineq_1-b}
\|\theta^{1}\|_{L^2(\Omega)} \leq  2\left( \bar{A}_{1} \tau^2 + \bar{B}_1 h^{s} + (1+\bar{C}_{1} \tau)\|\theta^{0}\|_{L^2(\Omega)} + C_{18}\tau \|\theta^{1^-}\|_{L^2(\Omega)} \right),
\end{equation}
and consequently
\begin{equation}\label{eq:ineq_1b}
\|\theta^{1}\|_{L^2(\Omega)} \leq  E_{1} ( \tau^2 + h^{s} )
\end{equation}
where $E_{1}$ does not depend on $\tau$ or $h$.

%%%%%%%%%%%%%
%%*********************
%%%%%%%%%%%%%

We proceed to analyze the case where $m\geq 2$. Let us define   
\begin{equation}\label{eq:def_Nset}
N_{h,\tau}^* = \max\left\{N^*\leq N={T}/{\tau}\,:\, N^*\in \mathcal{N}_{h,\tau}\right\}
\end{equation}
where $\mathcal{N}_{h,\tau} = \left\{ \bar{N} \,:\, \|{U}^k\|_{L^\infty(\Omega)} \leq 2\|v\|_{L^\infty([0,T]\times\Omega)}, \;k=0,\dots,\bar{N} \right\}$.  Observe that taking $h$ and $\tau$, sufficiently small, we have that $N_{h,\tau}^*\geq 1$. Indeed, we already proved \eqref{eq:U0inf_b}, and the bound $\|{U}^1\|_{L^\infty(\Omega)} \leq 2\|v\|_{L^\infty([0,T]\times\Omega)}$ can be established in the same way as \eqref{eq:U1-inf_b} by using \eqref{eq:ineq_1b}. 

Since $m\geq 2$, we have that $g^m = g(\widetilde{U}^m)$. Assuming that $m-1\leq N_{h,\tau}^*$,  we can easily get that $\|\widetilde{U}^m\|_{L^\infty(\Omega)} \leq 4 \|v\|_{L^\infty([0,T]\times\Omega)}$ and using \eqref{eq:ineq_g} it follows that  
\begin{equation}\label{eq:ineq_gc2}
|(g(v(t_{m-\frac{1}{2}})) - g^m ,\hat{\theta}^m)| \leq C_{21} \| v(t_{m-\frac{1}{2}}) - \widetilde{U}^m \|_{L^2(\Omega)} \|\hat{\theta}^m\|_{L^2(\Omega)}
\end{equation}
where $C_{21} = c_\rho (1+4^\rho)\|v\|_{L^\infty([0,T]\times\Omega)}$.

Consider the decomposition  
\[
v(t_{m-\frac{1}{2}}) - \widetilde{U}^m = \omega_5^m + \tilde{\eta}^m + \tilde{\theta}^m
\] 
where $\omega_5^m = v(t_{m-\frac{1}{2}}) - \tilde{v}^m$.
We have that
\[
\begin{aligned}
\omega_5^m &= \frac{3}{2}\left( \frac{\tau}{2} \partial_{t} v(t_{m-\frac{1}{2}})- \int_{t_{m-\frac{1}{2}}}^{t_{m-1}} (t_{m-1}-t) \partial_t^2 v(t)  dt\right) %\\&\qquad\qquad \qquad\qquad
- \frac{1}{2}\left(  \frac{3\tau}{2} \partial_{t} v(t_{m-\frac{1}{2}})- \int_{t_{m-\frac{1}{2}}}^{t_{m-2}} (t_{m-2}-t) \partial_t^2 v(t)  dt \right) \\
&= -\frac{3}{2} \int_{t_{m-\frac{1}{2}}}^{t_{m-1}} (t_{m-1}-t) \partial_t^2 v(t)  dt  %\\& \qquad\qquad 
+ \frac{1}{2} \int_{t_{m-\frac{1}{2}}}^{t_{m-2}} (t_{m-2}-t) \partial_t^2 v(t)  dt.
\end{aligned}
\]
Then,
\begin{equation}\label{eq:ineq_omega5}
%\begin{aligned}
\|\omega_5^m \|_{L^2(\Omega)} \leq  \frac{3\tau}{2}\left(\int_{t_{m-2}}^{t_{m-\frac{1}{2}}} \| \partial_t^2 v(t) \|_{L^2(\Omega)}  dt \right).
\end{equation}
We also have as a direct consequence of property \eqref{eq:est_proj_i} that
\begin{equation}\label{eq:ineq_eta2}
\begin{aligned}
\| \tilde{\eta}^{m}\|_{L^2(\Omega)} & \leq C_{22} h^s \left( \|v(t_{m-2})\|_{H^s(\Omega)} + \|v(t_{m-1})\|_{H^s(\Omega)} \right) \\ 
& \leq 2C_{22} h^s \|v\|_{L^\infty(0,T; H^s(\Omega))}.
\end{aligned}
\end{equation}
And the triangle inequality gives
\begin{equation}\label{eq:ineq_theta2}
\| \tilde{\theta}^{m}\|_{L^2(\Omega)} \leq  \frac{1}{2}\|\theta^{m-2}\|_{L^2(\Omega)} +\frac{3}{2} \|\theta^{m-1}\|_{L^2(\Omega)}.
\end{equation}

Finally, from equations \eqref{eq:ineq_gc2}--\eqref{eq:ineq_theta2}
 it follows that
 \begin{multline}\label{eq:S7}
 S^m_7  \leq 
 C_{23}  \Bigl( \tau \int_{t_{m-2}}^{t_{m-\frac{1}{2}}} \| \partial_t^2 v(t) \|^2_{L^2(\Omega)}  dt  + h^s \|v\|_{L^\infty(0,T; H^s(\Omega))}
 + \|\theta^{m-2}\|_{L^2(\Omega)} + \|\theta^{m-1}\|_{L^2(\Omega)}\Bigr) \\\times \left( \| \theta^m\|_{L^2(\Omega)} +  \| \theta^{m-1}\|_{L^2(\Omega)}  \right).
 \end{multline}

From equations \eqref{eq:S1}, \eqref{eq:S2}, \eqref{eq:S3}, \eqref{eq:S4}, \eqref{eq:S5}, \eqref{eq:S6} and \eqref{eq:S7}, we arrive to the conclusion that 
\begin{equation}\label{eq:ineq_ma}
\|\theta^m\|_{L^2(\Omega)} \leq  \bar{A}_m \tau^2 + \bar{B}_m h^{s} +C_{23} \tau \|\theta^{m-2}\|_{L^2(\Omega)} 
+ (1+\bar{C}_0\tau ) \| \theta^{m-1}\|_{L^2(\Omega)} +  \bar{D}\tau \| \theta^{m}\|_{L^2(\Omega)}
\end{equation}
where
\begin{equation}
\begin{split}
\bar{A}_m &=  C_{24} \int_{t_{m-2}}^{t_{m}} \bigl(\|  v(t)\|_{H^2(\Omega)}+ \|  \partial_{t} v(t)\|_{H^2(\Omega)}+ \|  \partial_t^2 v(t)\|_{H^2(\Omega)}  +  \| \partial_t^3 v(t)\|_{L^2(\Omega)} \bigr) dt, \\
\bar{B}_m &= C_{25}\Bigl( \tau \|v\|_{L^\infty(0,T;H^s(\Omega))}  + \int_{t_{m-1}}^{t_m}  \left( \|v(t) \|_{H^s(\Omega)}  +\|\partial_{t} v(t) \|_{H^s(\Omega)} \right) dt \Bigr), \\
\bar{C}_0 &= C_{23}+\bar{D}.
\end{split}
\end{equation}

From inequalities \eqref{eq:ineq_1a} and \eqref{eq:ineq_ma},  for  $1\leq m \leq \bar{N}_{h,\tau}$ with $ \bar{N}_{h,\tau} = \min\{N, N_{h,\tau}^* + 1\},$ we get that
\begin{equation}\label{eq:est-1_theta_m}
\|\theta^m\|_{L^2(\Omega)} \leq \bar{A} \tau^2 + \bar{B} h^s + \|\theta^0\|_{L^2(\Omega)}  + C_{18}\tau \|\theta^{1^-}\|_{L^2(\Omega)} 
+ \bar{C}\tau \sum_{j=0}^{m-1} \| \theta^{j}\|_{L^2(\Omega)} + \bar{D}\tau \| \theta^{m}\|^2_{L^2(\Omega)}
\end{equation}
where 
\[
\begin{aligned}
\bar{A} &= C_{26} \int_{0}^{T} \bigl(\|  v(t)\|_{H^2(\Omega)}+ \|  \partial_{t} v(t)\|_{H^2(\Omega)}+ \|  \partial_t^2 v(t)\|_{H^2(\Omega)}  +  \|\partial_t^3 v(t)\|_{L^2(\Omega)} \bigr) dt \\[8pt]
&\quad\qquad\leq C_{26} T^{1/2} \left( \|v\|_{W^2(0,T; H^2(\Omega))} + \|\partial_t^3 v\|_{L^2(0,T; L^2(\Omega))} \right),\\ 
\bar{B} &= C_{27}  \left(T \|v\|_{L^\infty(0,T;H^s(\Omega))} + \int_{0}^{T}  \left( \|v(t) \|_{H^s(\Omega)}  +\|\partial_{t} v(t) \|_{H^s(\Omega)} \right) dt \right)\\[8pt]
&\quad\qquad \leq C_{27} \left(2T \|v\|_{L^\infty(0,T;H^s(\Omega))} +  T^{1/2}\|\partial_{t} v\|_{L^2(0,T; H^s(\Omega))} \right), \\[8pt]
\bar{C} &= \max\{\bar{C}_1,\bar{C_0}\} + C_{23}.
\end{aligned}
\]
Notice that \eqref{eq:est-1_theta_m} is also valid for $m=0$. Moreover, by taking $\tau\leq (2\bar{D})^{-1}$ we obtain that
\begin{equation}\label{eq:est-2_theta_m}
\|\theta^m\|_{L^2(\Omega)} \leq 2\bigl( \bar{A} \tau^2 + \bar{B} h^s + \|\theta^0\|_{L^2(\Omega)}  + C_{18}\tau \|\theta^{1^-}\|_{L^2(\Omega)}  + \bar{C}\tau \sum_{j=0}^{m-1} \| \theta^{j}\|_{L^2(\Omega)}  \bigr)
\end{equation}
and the discrete Gronwall lemma \cite{Agarwal_book} implies that for $m=0, \dots,  \bar{N}_{h,\tau}$
\begin{equation*}%\label{eq:est-3_theta_m}
\|\theta^m\|_{L^2(\Omega)} \leq 2\left( \bar{A} \tau^2 + \bar{B} h^s + \|\theta^0\|_{L^2(\Omega)} +C_{18}\tau \|\theta^{1^-}\|_{L^2(\Omega)} \right) e^{2\bar{C}T}.
\end{equation*}
Hence, taking estimate \eqref{eq:ineq_1-b} into account, we conclude that
\begin{equation}\label{eq:est-3_theta_m}
\|\theta^m\|_{L^2(\Omega)} \leq \bar{E}( \tau^2 +  h^s) 
\end{equation}
where
\[
\bar{E} =  2e^{2\bar{C}T} \max\left\{ \bar{A} + \frac{C_{18} E_1^-}{2\bar{D}}, \bar{B} + C_0 + \frac{C_{18} E_1^-}{2\bar{D}} \right\}.
\]
Finally, applying property \eqref{eq:est_proj_i}, we bound $\|\eta^m\|_{L^2(\Omega)}$ and using the triangle inequality, we obtain the required estimates for $\|v(t_m)-U^m\|_{L^2(\Omega)}$ whenever $m=0,\dots,\bar{N}_{h,\tau}$.

In order to complete the proof, we shall establish that we can choose $h_0>0$ such that $\bar{N}_{h,\tau} = N = T/\tau$ whenever  $h\leq h_0$ and $\tau=o(h^{n/4})$. 
On the contrary, suppose this assertion is false, i.e. there is no such an $h_0$.  Then, we can always take a sufficiently small  $h$ and $\tau\leq (2\bar{D})^{-1}$ such that  $\bar{N}_{h,\tau} < N$. Therefore,  estimate \eqref{eq:est-3_theta_m} holds for $m = N^*_{h,\tau}+1 <N$. Moreover, since $h$  is sufficiently small and $\tau = o(h^{n/4})$, we can prove that $\|U^m\|_{L^\infty(\Omega)} \leq 2\|v\|_{L^\infty([0,T]\times\Omega)}$ in much the same way as estimate \eqref{eq:U1-inf_b}. This contradicts the definition of $N^*_{h,\tau}$, and the proof is complete.

\end{proof}

%%%%%%%%%%
\section{Numerical examples} \label{sec:num_ex}
This section will illustrate the error estimate result established in Theorem \ref{teo:caso_discreto} for the case where $s$ equals $r$.
Specifically, we will numerically demonstrate that the error between the exact solution to problem \eqref{prob:equivalente} and its approximate solution defined in \eqref{prob_aprox:discrete_all} satisfies the following error estimate:
\begin{align*}
    E(h,\tau):=
    \max_{m= 0,\dots, T/\tau} \|v(t_m)-U^m\|_{L^2(\Omega)}   
    \leq   C(\tau^2 + h^r).
\end{align*}
Additionally, we will present the approximate solution for various input data for the model under study.

We will illustrate convergence results and display the approximate solution for both one-dimensional ($1D$) and two-dimensional ($2D$) cases. 
In both dimensions, we will explore four distinct bases for the approximate subspace $S_h^r(\Omega)$: Lagrange polynomials of degrees 1, 2, and 3, and cubic Hermite polynomials.
It is worth noting that the relationship between the convergence order $\mathcal{O}(h^r)$ and the degree $p$ of the polynomials used as a basis is given by $r=p+1$, where $h$ represents the maximum diameter of the finite elements in the partition of $\Omega$.
For this study, we consider $\Omega = ]0,1[$ for the 1D case and $\Omega = ]0,1[\times ]0,1[$ for the 2D case, utilizing uniform partitions in both scenarios — meaning intervals of the same length in the 1D case and squares with equal sides in the 2D case.

\begin{figure}[H]
    \centering
    \includegraphics[width=0.9\textwidth]{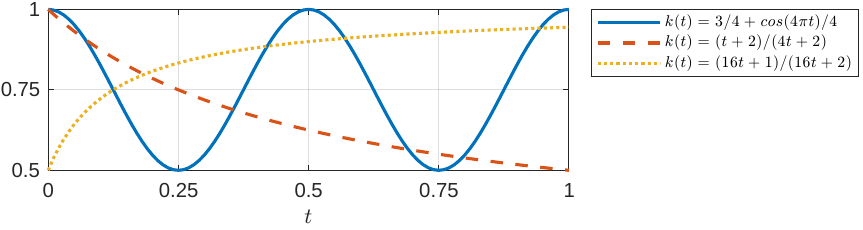}
    \caption{Functions $k(t)$ used to define the noncylindrical domains explored in this work.}
    \label{fig:Fig01}
\end{figure}

The noncylindrical domain at each time instant $t$ is represented by $\Omega_t\times{t}$, where $\Omega_{t}=\{ x\in\R^n;x=k(t)y, y \in \Omega\}$. 
The choice of $\Omega$ has already been established, leaving only the functions $k(t)$ specification to define the considered boundaries.
We will address three cases, which are defined below and visually represented in Figure \ref{fig:Fig01}.
\begin{equation}\label{eq:boundaries}
%\left\{
\begin{aligned}
\text{Boundary 1: }k(t) &= \frac{3}{4} + \frac{1}{4}\cos\big(4\pi t\big), \\%\quad
\text{Boundary 2: }k(t) &= \frac{t+2}{4t+2}, \\%\quad
\text{Boundary 3: }k(t) &= \frac{16t+1}{16t+2}.
\end{aligned}
%\right.
\end{equation}

We have set the following configurations for the remaining input parameters: $\rho=3$ in the nonlinear term $g(v)=|v|^\rho v$, a final time of $T=1$, and initial solutions $v_0$ and $U_0$ as specified below.
\begin{align}\label{eq:initial_cond}
& v_0(y) = 
\left\{
\begin{aligned}
& \sin(\pi y)(1+\iim),                 && 1D\text{ case},
\\
& \sin(\pi y_1) \sin(\pi y_2)(1+\iim), && 2D\text{ case},
\end{aligned}\right.
\qquad
\text{and}\quad U^0= P_h(0,v_0).
\end{align}

Regarding the definition of the source function $f(y,t)$ in equation \eqref{prob:equivalente}, since we are initially studying the error in the $1D$ and $2D$ cases, it is crucial to define this function so that the analytical solution to the model is known.
We will employ a common approach, where we define $f$ based on a sufficiently regular function $v$ that meets the initial and boundary conditions, 
i.e. a manufactured solution.
This allows us to determine $f$ by substituting the chosen $v$ into equation \eqref{prob:equivalente}.
Therefore, in the error study, we consider the function $f$ based on the following choices of $v$:
\begin{align*}
& v(y,t) = 
\left\{
\begin{aligned}
& \sin(\pi y)(1+\iim)\exp(-t),                 && 1D\text{ case},
\\
& \sin(\pi y_1) \sin(\pi y_2)(1+\iim)\exp(-t), && 2D\text{ case}.
\end{aligned}\right.
\end{align*}

\begin{figure}[H]
    \centering
    \includegraphics[width=0.8\textwidth]{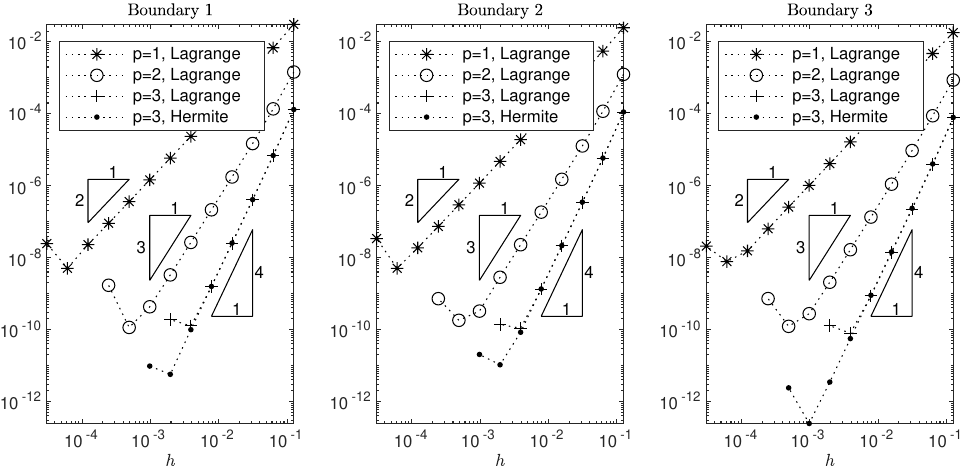}
    \caption{Non-homogeneous 1D problem.  Error plot $E(h,\tau)$ with $\tau=h^{\frac{p+1}{2}}$ for boundaries 1, 2 and 3.}
    \label{fig:Fig02}
\end{figure}

\begin{figure}[H]
    \centering
    \includegraphics[width=0.8\textwidth]{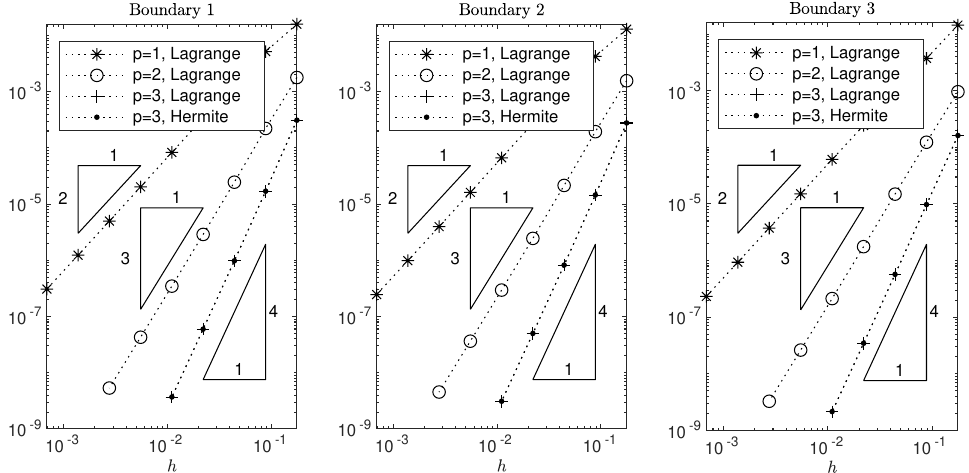}
    \caption{Non-homogeneous 2D problem. Error plot $E(h,\tau)$ with $\tau=h^{\frac{p+1}{2}}$ for boundaries 1, 2 and 3.}
    \label{fig:Fig03}
\end{figure}

In Figures \ref{fig:Fig02} and \ref{fig:Fig03}, we present the error analysis for the $1D$ and $2D$ cases, respectively.
Each curve shows the logarithmically scaled error $E(h,\tau)$ as a function of $h$, where $\tau=h^{(p+1)/2}$ and $p$ represents the degree of the polynomial basis of the finite element space.
By comparing the slope of the curves and the triangles, it is possible to conclude that the convergence order in space for each of the boundaries and the polinomial basis addressed is $\mathcal{O}(h^r)$, with $r$ being approximately $p+1$.
Regarding time convergence, since we are considering $\tau=h^{(p+1)/2}$, in the case of $p=1$, we have $\tau=h$. 
This means that for $p=1$, the error graph $E(h,\tau)$ can be analysed in terms of $h$ or $\tau$ since they are equivalent in this scenario.
As a result, from this case, we can conclude that the time convergence order is $\mathcal{O}(\tau^2)$. 
A similar conclusion applies to the other values of $p$. 


\begin{figure}[H]
    \centering
\includegraphics[width=\textwidth]{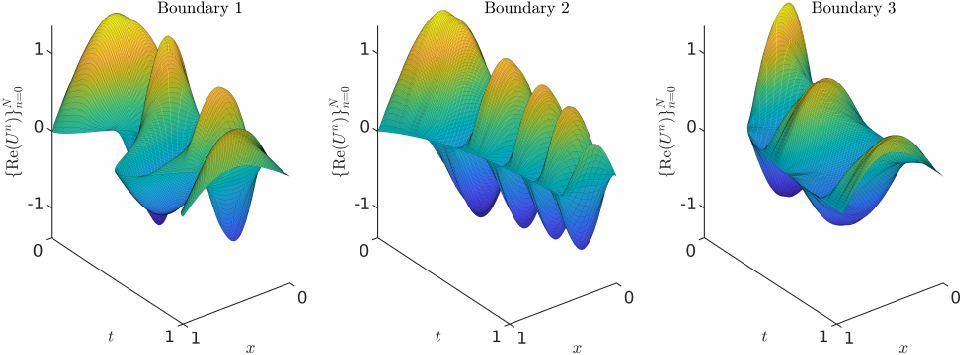}
\caption{Homogeneous 1D problem. Real part of approximate solutions relative to boundaries 1, 2, and 3, respectively. We consider cubic Hermite basis, $h=1/32$ and $\tau=1/400$.}
\label{fig:Fig04}
\end{figure}

\begin{figure}[H]
\centering
\includegraphics[width=\textwidth]{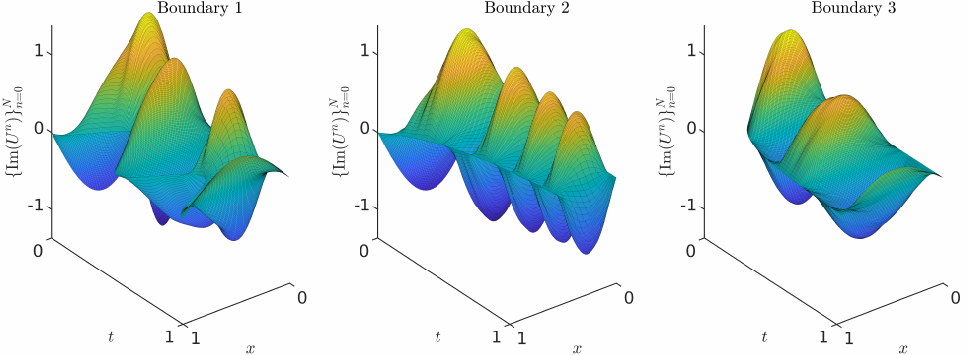}
\caption{Homogeneous 1D problem. Imaginary part of approximate solutions relative to boundaries 1, 2, and 3, respectively. We consider cubic Hermite basis, $h=1/32$ and $\tau=1/400$.}
\label{fig:Fig05}
\end{figure}

\begin{figure}[H]
    \centering
    \includegraphics[width=\textwidth]{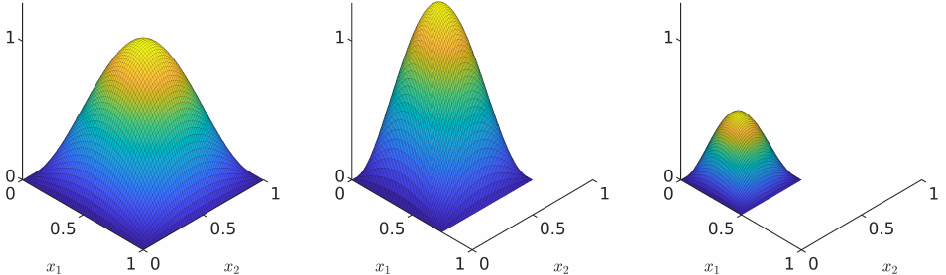}
    \caption{Homogeneous 2D problem. Real part of approximate solution $U^n$ relative to boundary 2 at times $t=0$, $t=0.25$, and $t=1$, respectively. We consider a cubic Hermite basis, $h=\sqrt{(1/32)^2+(1/32)^2}$ and $\tau=1/400$.}
    \label{fig:Fig06}
\end{figure}

\begin{figure}[H]
    \centering
    \includegraphics[width=\textwidth]{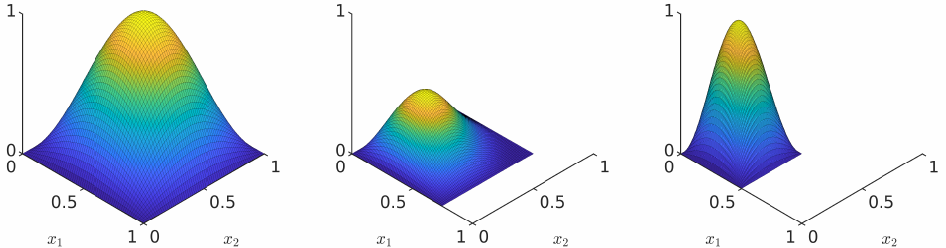}
    \caption{Homogeneous 2D problem. Imaginary part of approximate solution $U^n$ relative to boundary 2 at times $t=0$, $t=0.25$, and $t=1$, respectively. We consider a cubic Hermite basis, $h=\sqrt{(1/32)^2+(1/32)^2}$ and $\tau=1/400$.}
    \label{fig:Fig07}
\end{figure}

After concluding the error study, which provides us with confidence in the correctness of the approximate solution implementation, we now present in Figures \ref{fig:Fig04}-\ref{fig:Fig07} the approximate solution $U^m$ for the case where $f$ is identically zero, 
i.e. a homogeneous problem. The initial conditions are those used in the previous set of experiments and given by \eqref{eq:initial_cond}.
In this scenario, we no longer know the model's analytical solution. 

Figures \ref{fig:Fig04} and \ref{fig:Fig05} display the real and imaginary components of the approximate solution at each time step of the temporal discretization in the one-dimensional case. 
Analyzing these surfaces allows us to observe how the noncylindrical domain impacts the solution's behavior. 
These simulations were conducted using the cubic Hermite basis, with spatial discretization using $h=1/32$, and temporal discretization using $\tau=1/400$.

The following example corresponds to the bi-dimensional case with boundary 2. Figures \ref{fig:Fig06} and \ref{fig:Fig07} show the real and imaginary parts of the approximate solution, respectively, at different times. 
Unlike the one-dimensional case, where it is possible to display the solution over time discretization, allowing the visualization of the boundaries' influence on the solution, the same is not possible in the two-dimensional case. 
In fact, in the 2D case, the visualization at each time step \( t_n \), \( n=0,1,\dots,N \), is a surface. 
Therefore, in this work, we will display only the surfaces at three discrete times: \( t=0 \), \( t=0.25 \), and \( t=1 \). 
Despite this, as complementary and illustrative material, we provide the files accessible via the link \cite{code_gif}, along with videos of approximate solutions for the three exemplified boundaries.
    
\section{Conclusions}\label{sec:conc}

This article has mathematically established the convergence of approximate semidiscrete and fully discrete solutions to a nonlinear Schr\"odinger problem in a noncylindrical domain. For the fully discrete approximations given by a linearized Crank-Nicolson Galerkin method, we establish an optimal error estimate of order $\mathcal{O}( \tau^2 + h^{r})$, where $\tau$ and $h$ are the time step and the size of the spatial finite element mesh, respectively, and $r$ is the order of the finite element space.  Numerical simulations utilizing various basis demonstrate the consistency between theoretical predictions and numerical results for one-dimensional and two-dimensional cases. 

\section*{Acknowledgements}  
The author Mauro A. Rincon was partially supported by CNPq-Brazil and FAPERJ-Brazil (SEI 260003/003443/2022), and author Bruno A. do Carmo was supported by a post-doc grant from CAPES-Brazil (N. 23038.006308/2021-70).

\bibliographystyle{elsarticle-num}
\bibliography{bibliography}

\end{document}